\documentclass[a4paper]{article}

\usepackage{amsmath}
\usepackage{amsfonts}
\usepackage{amssymb}         %additional math symbols (see amsguide.tex)
\usepackage{amsopn}
\usepackage{theorem}          %additional LaTeX2e package (see pctex\latex2e)
\usepackage{latexsym}
\usepackage{graphicx}        %graphic and .eps files
\usepackage{mathrsfs}		%calligraphic fonts
\usepackage{psfrag}
\usepackage{algorithmic}
\usepackage{algorithm}
\usepackage{color}
\usepackage{verbatim}
\usepackage{enumerate}
\usepackage{url}

\theoremstyle{changebreak}                % (see LATEX2E\THEOREM.DTX)

 {{\nopagebreak\hspace*{\fill}$\Box$\par\vspace{12pt}}}
\begin{document}

\thispagestyle{empty}
\begin{center} 

{\LARGE Open research areas in distance geometry}
\par \bigskip
{\sc Leo Liberti${}^{1}$, Carlile Lavor${}^{2}$} 
\par \bigskip
\begin{minipage}{15cm}
\begin{flushleft}
{\small
\begin{itemize}
\item[${}^1$] {\it CNRS LIX, \'Ecole Polytechnique, F-91128 Palaiseau,
  France} \\ Email:\url{liberti@lix.polytechnique.fr}
\item[${}^2$] {\it IMECC, University of Campinas, Brazil} \\
  Email:\url{clavor@ime.unicamp.br}
\end{itemize}
}
\end{flushleft}
\end{minipage}
\par \medskip \today
\end{center}
\par \bigskip

% insert abstract
\begin{abstract} 
Distance Geometry is based on the inverse problem that asks to find the positions of points, in a Euclidean space of given dimension, that are compatible with a given set of distances. We briefly introduce the field, and discuss some open and promising research areas.
\end{abstract}

% insert paper

\section{Introduction}
\label{s:intro}
Distance Geometry (DG) is based around the concept of distance rather than points and lines. Its development as a branch of mathematics is mainly due to the motivating influence of other fields of science and technology, although pure mathematicians have also worked in DG over the years. DG becomes necessary whenever one can collect or estimate measurements for the pairwise distances between points in some set, and is then required to compute the coordinates of those points compatible with the distance measurements.

The fundamental problem in DG is the {\sc Distance Geometry Problem} (DGP), a decision problem that, given an integer $K>0$ and a connected simple edge-weighted graph $G=(V,E,d)$ where $d:E\to\mathbb{R}_+$, asks whether there exists a {\it realization} $x:V\to\mathbb{R}^K$ such that:
\begin{equation}
\forall \{u,v\}\in E \quad \|x_u-x_v\| = d_{uv},\label{eq:dgp}
\end{equation}
where $\|\cdot\|$ indicates an arbitrary norm, making this into a {\it problem schema} parametrized on the norm (as we shall see in the following, most applications employ the $\ell_2$ norm).

The DGP is an {\it inverse problem}: whereas computing some of the pairwise distances given the positions of the points is an easy task,\footnote{Direct problems may not be all that easy: see Erd\H{o}s' {\it unit distances} and {\it distinct distances} problems \cite{garibaldi}.} the inverse inference (retrieving the point positions given some of the distances) is not so easy. We remark that a realization can be represented by a $|V|\times K$ matrix, the $i$-th row of which is the position vector $x_i$ for vertex $i\in V$.

The main purpose of this paper is to survey what we think are the most important open problems in the field of DG.

The rest of this paper is organized as follows. We discuss some of the applications driving research behind DG in Sect.~\ref{s:apps}. We give a very short (and partial) historical overview of the development of DG in Sect.~\ref{s:history}. Sect.~\ref{s:probs} introduces some DGP variants. The main section, Sect.~\ref{s:open}, presents many open questions and promising areas for further research.

\section{Main application areas}
\label{s:apps}
The DGP arises in many application areas. We list here some of those for which the application is direct. In the following, we shall refer to $X$ as the set of solutions of the DGP variant being discussed.
\begin{itemize}
\item In certain network synchronization protocols, the time difference between certain clocks can be estimated and exchanged, but what is actually required is the absolute time \cite{singer4}. Here the points in the solution set $X$ are sequences of time instants, each of which is a vector in $\mathbb{R}^1$ (i.e.~a scalar) indicating the absolute time for a given clock. The time differences are Euclidean distances between points in one dimension: we recall that $\|\cdot\|_2=|\cdot|$ in $\mathbb{R}^1$. 
\item In wireless networks the devices may move about, usually on a two-dimensional surface. Some of the pairwise distances, typically those that are sufficiently close, can be estimated by measuring the battery power used in peer-to-peer communication. The information of interest to the network provider is the localization of the devices. In this setting, the solution set $X$ contains sequences of 2D coordinates (one per each device), and the measured distances are assumed to be approximately Euclidean, although they may be noisy and imprecise. See \cite{yemini78,yemini,doherty,jjsurvey,biswas2004,krislocksiam,eren04}
\item Proteins are organic molecules consisting of a backbone with some side-chains attached. Proteins have chemical properties (e.g.~the atoms which compose it) and geometrical properties (the relative position of each atom in the protein). Nowadays we know the chemical compositions of most proteins of interest, but not necessarily their shape; and yet, proteins bind to other proteins, and/or to specific sites on the surface of living cells, depending on both shape and chemical composition. Some of the pairwise inter-atomic distances can be measured in fairly precise ways (e.g.~the lengths of the covalent bonds). Kurt W\"uthrich discovered that some other distances (typically smaller than 5{\AA}) can be estimated using Nuclear Magnetic Resonance (NMR) techniques \cite{wuthrich_89}. In this setting, $X$ contains sequences of 3D coordinates (one per atom), and some of the measured distances can be noisy or just wrong \cite{berger}. Also see \cite{sidechains,nilges,bipbip}
\item Fleets of unmanned autonomous underwater vehicles (AUV) are deployed in order to install and maintain offshore installations such as oil rigs, wind energy farm and wave energy farms; such devices can estimate distances between each other, with the ocean bed and with the installation using sonar signals. In this setting, $X$ contains sequences of 3D coordinates (one per AUV). Measurements can be noisy and depend on time, and the positions of the AUVs move continuously in the ocean. See \cite{bahr}.
\item Nanostructures, such as graphite or buckminsterfullerene, are used extensively in material sciences. The main issue is determining their shape from a spectrographic analysis. The input data is essentially a probability distribution function (p.d.f.) over distance values, i.e.~a function $\mathbb{R}_+\to[0,1]$: by looking at the peaks of the p.d.f., one can extract a sequence of most likely values (with their multiplicities). From this list of distance values, one has to reconstruct their incidences to atoms, i.e.~the graph, and its realization in $\mathbb{R}^3$. In this setting, $X$ contains sequences of 3D coordinates (one per atom occurring in the nanostructure). The distances are {\it unassigned}, i.e.~they are simply values with multiplicities, but the incidence to the adjacent atoms is also unknown. See \cite{tribond,dg-4or,thorpe}.
\item In the analysis of robotic movements one is given the {\it bar-and-joint structure} of the robot (i.e.~a geometric model consisting of idealized rigid bars held together by freely rotating joints), the absolute position of its joints, and the coordinates of one or more points in space, and one would like to know if the robot can flex to reach that point. This involves computing the manifold of solutions of the DGP corresponding to the robotic graph weighted by the bar lengths, and asking whether the solution manifold contains the given points \cite{rojas_10}. Again, $X$ contains sequences of 3D coordinates (one per joint). 
\end{itemize}

In other cases, DG may be one of the steps towards the solution of other problems. In data analysis, for example, one often wishes to represent high-dimensional vectors visually on a 2D or 3D space. Many dimensional reduction techniques aim at approximately preserving the pairwise distances between the vectors \cite{borg_10}. This is fairly natural in that the ``shape'' of a set of points could be defined through the preservation of pairwise distances: two sets of points for which there exists a congruence mapping one to the other can definitely be stated to have ``the same shape''. Thus, an approximate congruence between two sets (such as the one defined in Eq.~\eqref{eq:jll}) might well be taken as a working definition of the two sets ``having approximately the same shape''. In dimensional reduction, the dimension $K$ of the target space is unspecified.

\section{Some historical notes about DG}
\label{s:history}
Arguably the first mathematical result that can be ascribed to DG is Heron's theorem for computing the area of a triangle given its side lengths \cite{heron}. This was further generalized by Arthur Cayley to simplices of arbitrary dimensions, the volume of which turns out to be equal a scaled determinant of a matrix which is a function of the side lengths of the simplex \cite{cayley1841}. Karl Menger proposed an axiomatization of DG \cite{menger28,menger31} that provided necessary and sufficient conditions for a metric space to be isometrically embeddable in a Euclidean space of finite dimension \cite{blumenthal61,bowers,havel}.

Much of the impact of DG on engineering applications will be discussed in the rest of the paper. In this section, we focus on two cases where the history of DG made an impact on modern mathematics. For more information, see \cite{six}.

\subsection{Impact of DG on rigidity}
Motivated by applications to statics and construction engineering, DG played a prominent role in the study of rigid structures, i.e.~bar-and-joint frameworks having congruences as their only continuous motions (a {\it bar-and-joint framework} is a bar-and-joint structure together with positions for the joints; or, equivalently, a pair $(G,x)$ of graph $G$ with an associated realization $x$).

Euler conjectured in 1766 \cite{euler1766} that all three-dimensional polyhedra are rigid. Cauchy provided a proof for strictly convex polyhedra \cite{cauchyrigid} (Cauchy's original proof contained two mistakes, subsequently corrected by Steinitz and Lebesgue), and Alexandrov \cite{alexandrov} extended the proof to all convex polyhedra. If polyhedra are defined by their face incidence lattice rather than as intersections of half-spaces, then polyhedra can also be nonconvex: this, in fact, appeared to be the setting proposed by Euler in his original conjecture, expressed in terms of (triangulated) surfaces. In this setting, Bob Connelly finally found in 1978 \cite{connelly-countereg} an example of a nonconvex triangulated sphere which can undergo a flexible motion of some of its vertices, whilst keeping all the edge distances equal, and disproved Euler's conjecture.

J.C.~Maxwell studied rigidity \cite{maxwell1864,maxwell1864b} in relation to balancing forces acting on structures, more precisely force diagrams by reciprocal figures. These were at the basis of graphical algorithms to verify force balancing \cite{cremona1872,cremona1874}, in use until computers became dominant \cite{recskisurvey1}.

\subsection{The role of DG in ``big data''}
\label{s:bigdata}
DG is also at the center of two results currently used in the analysis and algorithmics of large data sets, also known as ``big data''. Both results gave rise to dimensional reduction techniques, i.e.~methods for projecting a finite set $Y$ of points in $\mathbb{R}^m$ (with $m$ large) to $\mathbb{R}^K$ (with $K$ much smaller than $m$), while approximately preserving the pairwise distances over $Y$. The first such technique is Multidimensional Scaling (MDS), originally based on a 1935 result of Isaac Schoenberg \cite{schoenberg} (rediscovered in 1938 by Young and Householder \cite{young_householder}). The second technique is based on a lemma of Johnson and Lindenstrauss \cite{jllemma}, which ensures that, for a given $\varepsilon\in(0,1)$, a large enough $|Y|$, and $K=O(\frac{1}{\varepsilon^2}\ln |Y|)$, there exists a function $f:\mathbb{R}^m\to\mathbb{R}^K$ which preserves pairwise distances approximately, up to a multiplicative error:
\begin{equation}
  \forall x,y\in Y\quad (1-\varepsilon)\|x-y\|_2 \le \|f(x)-f(y)\|_2 \le (1+\varepsilon) \|x-y\|_2.\label{eq:jll}
\end{equation}
MDS is now a pervasive data analysis technique, applied to a vast range of problems from science and technology. The Johnson-Lindenstrauss lemma (JLL) is less well known, but employed e.g.~for fast clustering of Euclidean data \cite{indyk}. A possible application of these results to the DGP is useful to project large dimensional realizations to a smaller dimension while keeping all pairwise distances approximately equal.

\section{Problems in DG}
\label{s:probs}
As already mentioned, the DGP is the fundamental problem in DG. In the DGP formulation we gave in Sect.~\ref{s:intro}, however, we omitted to specify the norm, which is mainly intended to be Euclidean. In this case, the DGP is also called {\sc Euclidean DGP} (EDGP) \cite{dgp-sirev,vetterli}.

In this section we look at several types of DGP variants. In Sect.~\ref{s:fixeddim} we look at the case of fixed $K$ given as part of the input. In Sect.~\ref{s:dgpnorms} we discuss the DGP using other norms than the Euclidean one. In Sect.~\ref{s:idgp} we discuss the DGP in the presence of measurement errors on the input data. In Sect.~\ref{s:isomemb} we discuss the case where $G$ is complete and $K$ is not given as part of the input, but rather as an asymptotic bound in function of $n=|V|$. In Sect.~\ref{s:matcompl} we discuss the case where $K$ is part of the output. Finally, in Sect.~\ref{s:unassigned} we present the DGP variant where the weighted graph is replaced by a list of distance values.

\subsection{DGP in given dimensions}
\label{s:fixeddim}
Although the dimension $K$ is specified in the DGP as part of the input, most applications require a fixed given constant, see Sect.~\ref{s:apps}: for example, the determination of protein structure from distance data requires $K=3$. When the dimension $K$ is fixed to a constant $\gamma$, we denote the corresponding problem by DGP${}_{\gamma}$ (equivalently, we denote by EDGP${}_\gamma$ the EDGP in fixed dimension $\gamma$, and similarly for other DGP variants). Because of the application to molecules, the EDGP${}_3$ is also called {\sc Molecular DGP} (MDGP); similarly, because of the application to wireless networks, the DGP${}_2$ is also called the {\sc Sensor Network Localization Problem} (SNLP). 

\subsection{DGP with different norms}
\label{s:dgpnorms}
The DGP using other norms is not as well studied as the EDGP. Some Mixed-Integer Linear Programming (MILP) formulations and some heuristics are currently being developed for the $\ell_1$ and $\ell_\infty$ norms \cite{oneinfnorm}. The $\ell_\infty$ norm is used as a proxy to the $\ell_2$ norm in \cite{crippeninf}. Some works in spatial logic reduce to solving the DGP with a discrete-valued semimetric taking values in a set such as $\{\mbox{\tt almost\_equal},\mbox{\tt near},\mbox{\tt far}\}$ (each label is mapped to an interval of possible Euclidean distances) \cite{splogic}. Geodesic spherical distances have been briefly investigated by G\"odel \cite{goedelDG1}, who proved that if a weighted complete graph over 4 vertices can be realized in $\mathbb{R}^3$, but not in $\mathbb{R}^2$, then it can also be realized on the surface of a sphere. This was extended to arbitrary dimensions in \cite{dgpsphere}.

\subsection{DGP with intervals}
\label{s:idgp}
In most applications, distances are not given precisely. Instead, some kind of measurement errors are involved. A common way to deal with this issue is to model distances $d_{uv}$ by means of an uncertainty interval $[d_{uv}^L,d_{uv}^U]$, yielding what is known as {\it interval} DGP ({\it i}DGP):
\begin{equation}
  \forall \{u,v\}\in E\quad d_{uv}^L\le \|x_u-x_v\|\le d_{uv}^U. \label{eq:idgp}
\end{equation}
Very few combinatorial techniques for solving DGPs extend naturally in the case of intervals. Typically, optimization techniques based on Mathematical Programming (MP) do, however \cite{morewu,dgpzoo-tr}. See \cite{mdgpsurvey,bpinterval,souza2,bipbip} for more information. 

\subsection{Isometric embeddings}
\label{s:isomemb}
Many works exist in the literature about {\it isometric embeddings}, i.e.~embeddings of metrics in various vector spaces. We look specifically at cases where the metrics are finite and the target space is the Euclidean space $\mathbb{R}^K$ (for some $K$). We remark that the isometric embedding problem with finite metrics is close to the case of the DGP where the input graph $G$ is a complete graph.

In this line of (mostly theoretical) research, $K$ is not usually given as part of the input but rather proved to be a function of $|V|$ which is asymptotically bounded above. The JLL (Sect.~\ref{s:bigdata}) is an example of this type of results. An ingenious construction shows that any valid metric $D=(d_{uv})$ can be embedded exactly in $\mathbb{R}^n$ (where $n=|V|$) under the $\ell_\infty$ norm. It suffices to define \cite{frechet,kuratowski}:
\begin{equation}
  \forall v\in V \quad x_v = (d_{uv}\;|\;u\in V).\label{eq:frechet}
\end{equation}
This construction is known as the {\it Fr\'echet embedding}.

For the $\ell_1$ norm, no such general result is known. It is known that $\ell_2$ metric spaces of $n$ points can be embedded in a vector space of $O(n)$ dimensions under the $\ell_1$ norm almost isometrically \cite[\S 2.5]{matousekmetric}. The ``almost'' in the previous sentence refers to a multiplicative distortion similar to the JLL's:
\[(1-\varepsilon)\|x\|_2 \le \|f(x)\|_1\le (1+\varepsilon)\|x\|_2,\]
where $f$ preserves norms approximately while reducing the dimensionality of $x$, for some $\varepsilon\in(0,1)$. Moreover, any finite metric on $n$ points can be embedded in an exponentially large dimensional vector space using the $\ell_1$ norm with $O(\log n)$ distortion: this was shown in \cite{bourgain} by means of a probabilistic weighted extension of the Fr\'echet embedding on all subsets of $V$. The dimension was reduced to $O(n^2)$ by a deterministic construction \cite{linial}; moreover, an appropriate randomized choice of subsets of $V$ drives it down to $O(\log n)$ \cite{linial-stoc}. Similar results hold for many other norms, including $\ell_2$ one. The relatively large distortion $O(\log n)$ unfortunately limits the usefulness of these results.

\subsection{Matrix completion}
\label{s:matcompl}
The EDGP can also be formulated as follows: given $K>0$ and the squared weighted adjacency matrix $D=(d_{uv}^2)$ of $G$ (with $d_{uv}$ being the weight of the edge $\{u,v\}$ if $\{u,v\}\in E$), find a squared Euclidean Distance Matrix (sqEDM) $\bar{D}=(\bar{d}_{uv}^2)$ corresponding to a realization in $\mathbb{R}^K$ and such that, for each $\{u,v\}\in E$, $\bar{d}_{uv}=d_{uv}$.

The {\sc EDM Completion Problem} (EDMCP) consists in relaxing the requirement that $\bar{D}$ should be the sqEDM of a realization in $\mathbb{R}^K$ for a given $K$. Instead, $K$ is not part of the input, and $\bar{D}$ should simply correspond to a realization in a Euclidean space of any dimension. Informally, this means that, in the EDMCP, $K$ becomes (implicitly) part of the output. More details can be found in \cite{mcp}.

\subsection{DGP without adjacency information}
\label{s:unassigned}
Suppose that, instead of providing a weighted graph $(G,d)$, we provided instead a list $L$ of distance {\it values}, and then asked the same question. We can no longer write Eq.~\eqref{eq:dgp} since we do not know what $d_{uv}$ is: instead, we only have $L=(d_1,\ldots,d_m)$ where $m=|E|$. This DGP variant, called the {\sc unassigned DGP} (uDGP) is very important, since NMR and X-ray crystallography experiments actually provide the distance values rather than the actual edges of the graph. Although a fair amount of work has been carried out by physicists \cite{tribond} and structural biologists \cite{nilges3} on this problem, it is largely unstudied by mathematicians and computer scientists for all $K>1$ (the case $K=1$, on the other hand, has been studied under different names: turnpike problem and partial digest problem \cite{skiena}). This prompted us to list it as one of the main ``open areas'' in Sect.~\ref{s:open} below (see Sect.~\ref{s:udgp}).

\section{Open research areas}
\label{s:open}

In the last ten years, our work in the DG research community allowed us to survey many theoretical, methodological and applicative areas connected to DG \cite{lln4,mdgpsurvey,dmdgpejor,dgpbook,dgp-sirev,dgpdampreface,dgpitorpreface,dgta16proc}. Although our viewpoint is certainly not exhaustive, we list here the research topics which we think are most promising for further research.
\begin{enumerate}
\item Combinatorial characterization of rigidity in dimensions $K>2$.
\item Computational complexity of Euclidean Distance Matrix Completion in the Turing Machine model.
\item {\it A priori} estimation of the number of realizations for a given set of distances.
\item The unassigned DGP.
\end{enumerate}

We also remark that the longest-standing problem in DG, that of the rigidity of all closed triangulated surfaces, was opened by L.~Euler in 1766 \cite{euler1766} and finally answered in the negative by Bob Connelly in 1978 \cite{connelly-countereg}.  

In the rest of this paper, we discuss each of these research areas in detail.

\subsection{Combinatorial characterization of rigidity}
\label{s:rigidity}
Rigidity is a property relating to the bar-and-joint framework theory of statics. Architects and construction engineers are concerned with structures that do not bend: or, in other words, that are rigid. From the point of view of many other applications employing the EDGP as a model of the inverse problem of recovering positions from distance information, a desirable property is solution uniqueness. For example, if one is trying to recover the position of wireless devices in a mobile network, one would like the solution to be unique. In the case of protein conformation, one would like to find all of the possible chiral isomers, which are in finite number. Again, the property that tells apart EDGP instances with a finite number of solutions from those with an uncountable number is rigidity.

Formally speaking, there are many different definitions of rigidity. The most basic one is concerned with lack of local movement, i.e.~the only possible movements that a structure can undergo without changing the given pairwise distances are rotations and translations. Another definition, most often used in statics, concerns the absence of infinitesimal motions (defined below). Other definitions concern solution uniqueness (global/universal rigidity) \cite{alfakihdgp}, minimality of edge set cardinality (isostaticity/minimal rigidity) \cite{tay-whiteley}, abstractness (graphical/abstract rigidity matroids) \cite{graverbook}, and more \cite{sitharam}.

If the framework has certain genericity properties, then rigidity can be ascribed directly to the graph, rather than the framework \cite{gluck}. The question then becomes: given a graph, is it rigid in $K$ dimensions? Since the input only consists of a graph and an integer, the ideal solution algorithm to settle this question should be ``purely combinatorial'', meaning that during its execution on a Turing Machine computational model, no real number should be approximated through floating point representations. Purely combinatorial characterizations are known for $K\in\{1,2\}$, but not for any other value of $K$. Currently, this is considered the most important open question in rigidity, and possibly for the whole of DG.

\subsubsection{Rigidity of frameworks}
\label{s:rigidfmk}
Consider a YES instance of the EDGP, consisting of a weighted graph $G=(V,E,d)$ and an integer $K>0$, as well as a realization $x\in\mathbb{R}^{Kn}$. The pair $(G,x)$ is called a {\it framework} in $\mathbb{R}^K$. We let $\mathbf{K}(G,x)$ be the complete graph over $V$ weighted by the edge function $\bar{d}$ defined as follows: $\bar{d}_{uv}=d_{uv}$ for all $\{u,v\}\in E$, and $\bar{d}_{uv}=\|x_u-x_v\|_2$ otherwise (we shorten $\mathbf{K}(G,x)$ to $\mathbf{K}$ when no ambiguities arise). We further define the {\it edge weight value function} $f_G:\mathbb{R}^{Kn}\to\mathbb{R}^{|E|}_+$ by $f_G(x) = (\|x_u-x_v\|_2\;|\;\{u,v\}\in E)$.

A framework $(G,x)$ is {\it rigid} if there is at least a neighbourhood $\chi$ of $x$ in $\mathbb{R}^{Kn}$ such that:
\begin{equation}
  f^{-1}_G(f_G(x))\cap\chi = f^{-1}_{\mathbf{K}}(f_{\mathbf{K}}(x)). \label{eq:rigidity}
\end{equation}
The expression $f^{-1}_G(f_G(x))$ denotes the set of realizations $x'$ that satisfy the same distance equations as $x$, i.e.~the set of all solutions of the given EDGP instance. The LHS therefore indicates all those realizations of the EDGP instance which are in the neighbourhood of $x$. Similarly, the RHS indicates the same when $G$ is replaced by its {\it completion} $\mathbf{K}(G,x)$.

Realizations of complete graphs can be moved isometrically (i.e.~while keeping the edge lengths invariant) only if the movements are {\it congruences}, namely compositions of rotations and translations (we do not consider reflections since they are non-local). The intuitive sense of the above definition is that, if Eq.~\eqref{eq:rigidity} holds, then $G$ locally ``behaves like'' its completion. In other words, a framework is rigid if it can only be moved isometrically by congruences. Testing rigidity of a given framework (with a rational realization) when $K=2$ is co{\bf NP}-hard \cite{abbott} (see Sect.~\ref{s:complbasic} for a definition of co{\bf NP}). 

\subsubsection{Infinitesimal rigidity}
\label{s:infrigid}
We now focus on rigidity from the point of view of the movement. Consider the graph \[G=(\{1,2,3\},\{\{1,2\},\{2,3\}\})\] on three vertices, with two edges: node 2 is adjacent to both 1 and 3, and hence has {\it degree} 2, while nodes 1 and 3 both have degree 1. Consider the realization $x_1=(0,0)$, $x_2=(1,0)$, and $x_3=(1,1)$ in $\mathbb{R}^2$. It is obvious that any position for $x_3$ on the unit circle centered at $x_2$ is also a valid realization: this generates an uncountable set $x(\alpha)$ of realizations as the angle $\alpha$ between the segments $\overline{12}$ and $\overline{23}$ ranges in the interval $[0,2\pi]$. If one sees $\tau$ as a time parameter, the variation of $x_3$ is an isometric movement, also known as a (nontrivial) {\it flex}, which implies that this graph is flexible for $K=2$ (by contrast, it is rigid for $K=1$). 

We now generalize this to any graph $G=(V,E)$ with any realization $x$ in $\mathbb{R}^K$. By isometry, any flex has to satisfy Eq.~\eqref{eq:dgp}, which we write in its squared form:
\begin{equation*}
  \forall \tau\in [0,1], \{u,v\}\in E \quad \|x_u(\tau)-x_v(\tau)\|_2^2 = d_{uv}^2. 
\end{equation*}
Note that we can assume that $\tau\in[0,1]$ by rescaling if necessary. Since $x_u(\tau)$ denotes the position in $\mathbb{R}^K$ of vertex $u$ at time $\tau$, we can compute its velocity by simply taking derivatives with respect to $\tau$:
\begin{equation*}
  \forall \tau\in [0,1], \{u,v\}\in E \quad \frac{\mbox{\sf d}}{\mbox{\sf d}\tau} \|x_u(\tau)-x_v(\tau)\|_2^2 = \frac{\mbox{\sf d}}{\mbox{\sf d}\tau} d_{uv}^2.
\end{equation*}
We now remark that the RHS of this equation is zero, since $d_{uv}^2$ are constants for all $\{u,v\}\in E$, yielding:
\begin{equation*}
  \forall \{u,v\}\in E \quad  (x_u-x_v)\cdot(\dot{x}_u-\dot{x}_v) = 0,
\end{equation*}
where we assume that $x=x(0)$ is the realization given in the framework $(G,x)$.
In order to find the velocity vector $\dot{x}$ at $\tau=0$, we have to compute $R_{uv}=x_u-x_v$ for all $\{u,v\}\in E$, then solve the homogeneous linear system
\begin{equation}
  R\, \dot{x} = 0, \label{eq:infrigid}
\end{equation}
where $R$ is a matrix having $|E|$ rows and $Kn$ columns, called {\it rigidity matrix}. $R$ is defined as follows: for a row indexed by $\{u,v\}\in E$, there are $K$ possibly nonzero columns indexed by $(u,1),\ldots,(u,K)$ containing the entries
\[ x_{u1} - x_{v1},\ldots, x_{uK} - x_{vK},\]
and $K$ possibly nonzero columns indexed by $(v,1),\ldots,(v,K)$ containing the reciprocal entries
\[ x_{v1} - x_{u1},\ldots,x_{vK}-x_{uK}.\]
Sometimes rigidity matrices are shown in ``short-hand format'' $|E|\times n$ by writing each sequence of entries $(x_{u1}-x_{v1},\ldots,x_{uK}-x_{vK})$ as $x_u-x_v$ in the column indexed by $u$, and equivalently as $x_v-x_u$ in the column indexed by $v$. 

We now make the following crucial observation: the vector subspace spanned by all $\dot{x}$ satisfying Eq.~\eqref{eq:infrigid} contains all of the instantaneous velocity vectors that yield isometric movements of the given framework, also called {\it infinitesimal motions} of the framework. We remark that this subspace corresponds to the kernel $\mbox{\sf ker}\,R$ of the rigidity matrix. Now, the framework $(G,x)$ is {\it infinitesimally rigid} if $\mbox{\sf ker}\,R$ only encodes the translations and rotations of $(G,x)$. Otherwise the framework is {\it infinitesimally flexible}.

Since infinitesimal rigidity is based on a linear system, it can be decided based on the estimation of the degrees of freedom of the framework. If we start with an empty edge set, each of the $n$ vertices has $K$ degrees of freedom, so the system has $Kn$ degrees of freedom overall. Each linearly independent row of $R$ decreases the degrees of freedom by one unit. In general, the framework $(G,x)$ has $Kn-\mbox{\sf rk}\,R$ degrees of freedom (we denote the rank of $R$ by $\mbox{\sf rk}\,R$). We remark that in $\mathbb{R}^K$ there are $K$ basic translations and $K(K-1)/2$ basic rotations (arising from pairs of distinct orthogonal axes), for a total of $K(K+1)/2$ degrees of freedom of the group of Euclidean motions in $\mathbb{R}^K$. Since any framework can be moved isometrically by this group, at least these basic motions must satisfy Eq.~\eqref{eq:infrigid}. Hence, we have
\[\mbox{\sf dim\,ker}\,R\ge K(K+1)/2.\]
Specifically, $(G,x)$ is infinitesimally rigid if and only if
\begin{equation}
  \mbox{\sf rk}\,R = Kn - \frac{K(K+1)}{2} \label{eq:asimowA}
\end{equation}
and infinitesimally flexible if and only if
\begin{equation}
  \mbox{\sf rk}\,R < Kn - \frac{K(K+1)}{2},\label{eq:asimowB}
\end{equation}
as long as the affine hull of $x$ spans the whole of $\mathbb{R}^K$ \cite{asimow1}. Moreover, it was shown in \cite{asimow2} that a framework is infinitesimally rigid if and only if it is rigid and $x$ is {\it regular}, meaning that its rigidity matrix has the maximum possible rank.

We remark that infinitesimal rigidity is a stronger notion than rigidity: every infinitesimally rigid framework is also rigid, but there are rigid frameworks that are infinitesimally flexible. 

\subsubsection{Generic properties}
\label{s:genericprop}
Another consequence of Eq.~\eqref{eq:asimowA}-\eqref{eq:asimowB} is that, if one can find a single rigid realization of the graph $G$, then almost all realizations of $G$ must be rigid. This follows because: (a) rigidity is the same as infinitesimal rigidity with $R$ having maximum rank (say) $r$ among all the rigidity matrices corresponding to rigid frameworks; and (b) if a random matrix with a certain sparsity pattern is sampled uniformly from a compact set, it ends up having its maximum possible rank with probability 1. The alternative, i.e.~that a sampled $R$ is rank deficient, corresponds to the existence of a linear relationship between its rows, which defines a subspace of zero Lebesgue measure in $\mathbb{R}^r$.

Because of this, both rigidity and infinitesimal rigidity are {\it generic properties}. This allows us to ascribe them directly to the underlying (unweighted) graph rather than the framework. Equivalently, it is sufficient to look at the sparsity structure of the rigidity matrix in order to decide whether a framework is rigid with probability one.

Formally stated, we are concerned with the following decision problem.
\begin{quote}
  {\sc Graph Rigidity}. Given a graph $G=(V,E)$ and an integer $K>0$, is $G$ generically (infinitesimally) rigid in $K$ dimensions?
\end{quote}

It should be clear that we can potentially use the rank of the rigidity matrix in order to solve the problem. Consider this algorithm:
\begin{enumerate}
\item sample a random realization $x$ of $G$
\item compute the rank of $R$
\item if it is equal to $Kn-K(K+1)/2$ then output YES
\item otherwise output NO.
\end{enumerate}
This algorithm runs in polynomial time, since matrix ranks can be computed in polynomial time \cite{cheung}. But this is not considered a ``combinatorial characterization'' since it involves computation with floating point numbers.

Moreover, it is a randomized algorithm in the Turing Machine (TM) computational model, since only finitely many rationals can be represented in the commonly used IEEE 754 standard floating point implementation: sampling a random realization $x$ will yield a maximum rank rigidity matrix with practically high probability, but not with probability 1. Even with theoretical probability 1 assumptions, probability 1 is not the same as certainty: there is a possibility that this algorithm might output NO on a small fraction of YES instances. Re-running the algorithm sufficiently many times on the same instance will make the probability of error as small as desired, but this is not a deterministic algorithm. 

An acceptable ``combinatorial characterization'' would limit its scope to decision algorithms that only employ the incidence structure of $G$, and perhaps integers bounded by a polynomial in $|V|$ and $|E|$. Although the original meaning ascribed to ``combinatorial characterization'' did not call for algorithmic efficiency (in the sense of worst-case polynomial running time), the rigidity research community appears to feel that this would be a desirable characteristic \cite{dgta16proc}.

\subsubsection{Combinatorial characterization of rigidity on the line}
\label{s:combcharact1}
For $K=1$, a combinatorial characterization of generic rigidity is readily available: $G$ is generically rigid if and only if it is connected. This holds because the only flexes in $\mathbb{R}^1$ are translations. If a graph has a flex, say on vertex $v$, then all of its neighbours must undergo the same flex. If the graph is connected, then an induction argument shows that all of the vertices must undergo the same flex, showing that the flex is a congruence. If the graph is disconnected, then two connected components can undergo different flexes, which means different translations, the combination of which is not a congruence of the whole graph.

We remark that graph connectedness can be decided in polynomial time, for example by graph exploration algorithms \cite{mehlhorn}.

\subsubsection{Combinatorial characterization of rigidity in the plane}
\label{s:combcharact2}
For $K=2$, the situation becomes much more complicated. James Clerk Maxwell was using a degree of freedom calculus already in 1864 \cite{maxwell1864b} to the effect that if $G$ is minimally rigid in the plane, then $|E|=2|V|-3$ must hold. We remark that a graph is {\it minimally rigid} (also known as {\it isostatic}) when it is no longer rigid after the removal of any of its edges.

Gerard Laman finally proved his celebrated theorem in 1970 \cite{laman}, namely that $G=(V,E)$ with $|V|>1$ is generically minimally rigid if and only if:
\begin{itemize}
\item[(a)] $|E|=2|V|-3$
\item[(b)] for each subgraph $G'=(V',E')$ of $G$ having $|V'|>1$, we have $|E'|\le 2|V'|-3$.
\end{itemize}
This was accepted as a purely combinatorial characterization of rigidity in the plane, since it immediately gives rise to the following brute force algorithm. Given any graph $G$,
\begin{enumerate}
  \item list all subgraphs of $G$ with at least two vertices
  \item for each of them, test whether Laman's conditions (a) and (b) hold
  \item if they do, output YES, otherwise NO.
\end{enumerate}
Since there are exponentially many subgraphs of any graph, this algorithm is exponential time in the worst case. On the other hand, a polynomial time algorithm based on Laman's theorem was given in \cite{lovasz-yemini}.

\subsubsection{Combinatorial characterization of rigidity in space}
\label{s:combcharact3}
\begin{figure}[!ht]
  \begin{center}
    \includegraphics[width=9cm]{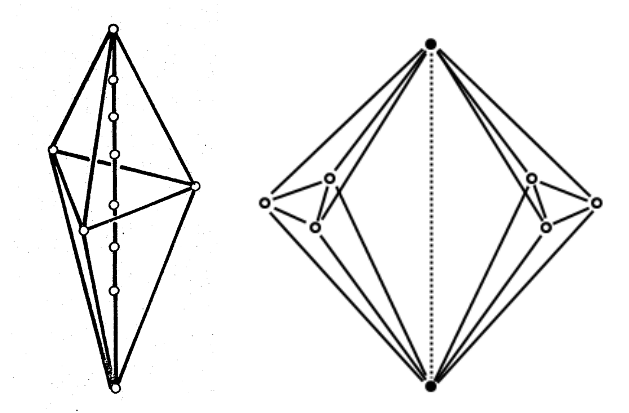}
  \end{center}
  \caption{Counterexamples for Laman's conditions in 3D. On the left, a rigid graph with fewer than $3|V|-6$ edges. On the right, a flexible graph satisfying the 3D adaptation of Laman's condition (b): the dotted line shows a hinge around which the two ``bananas'' can rotate.}
  \label{f:counterex}
\end{figure}

Although Maxwell did use the rank formula in 3D, namely $|E|=3|V|-6$, in one of his papers \cite{maxwell1864}, all attempts to extend Laman's theorem along these lines for the case $K=3$ have failed so far. In fact, Laman's conditions fail spectacularly when $K=3$. The equivalent of condition (a) above would be $|E|=3|V|-6$, but a well-contrived example \cite[p.~2.14]{grunbaum} (see Fig.~\ref{f:counterex}, left) shows this to be false. The failure of the equivalent of condition (b) above, i.e.~that every subgraph with $|V'|>1$ should have $|E'|\le 3|V'|-6$ is exhibited by the famous ``double banana'' graph (see Fig.~\ref{f:counterex}, right), ascribed to Whiteley by \cite{asimow2}. 

Two ideas which might shed some light over the long-standing open question of finding a combinatorial characterization of 3D (generic) rigidity have been proposed by Meera Sitharam. In 2004, she and Y.~Zhou have published a paper \cite{sitharam} based on the concept of module rigidity (to replace Laman-style degree of freedom counts), according to which graphs like the double banana of Fig.~\ref{f:counterex} (right) would not be erroneously catalogued as rigid while being flexible. More recently, her talk (co-authored with J.~Cheng and A.~Vince) at the Geometric Rigidity 2016 workshop in Edinburgh bore the title ``Refutation of a maximality conjecture or a combinatorial characterization of 3D rigidity?''. We report the abstract from the workshop proceedings, since there appear to be no publications about this idea yet.
\begin{quote}
  The talk will present an explicit, purely combinatorial algorithm that defines closure in an abstract 3D rigidity matroid (GEM, for graded exchange matroid). Strangely, rank in this matroid is an upper bound on rank in the 3D rigidity matroid, either refuting a well-known maximality conjecture about abstract rigidity matroids, (more likely) or giving a purely combinatorial characterization of 3D rigidity (less likely). In addition, we can show that rank of a graph G in the new GEM matroid is upper bounded by the size of any maximal $(3,6)$-sparse subset.
\end{quote}
Although she does warn that the most likely possibility is that her result is {\it not} the sought-after characterization, we believe that explicitly and visibly displaying attempts to solve hard problems, such as this one, has the advantage of drawing attention (and hence further study and effort) towards their solution. 

\subsubsection{Global rigidity}
\label{s:globrig}
A framework $(G,x)$ is {\it globally rigid} in $\mathbb{R}^K$ if it is rigid in $\mathbb{R}^K$ and only has one possible realization up to congruences (including reflections). This case is very important in view of applications to clock synchronization protocols, wireless network realization, autonomous underwater vehicles, and, in general, whenever one tries to estimate a unique localization occurring in the real physical world.

By \cite{gortler}, global rigidity is a generic property, i.e.~it depends only on the underlying graph $G$ (rather than the framework) for almost all realizations $x$. For $K=1$, $G$ is (generically) globally rigid if and only if it is bi-connected (i.e.~there are at least two distinct simple paths joining each pair of vertices in the graph). For $K=2$, $G$ is (generically) globally rigid if and only if it is 3-connected and redundantly rigid (i.e.~it remains rigid after the removal of any edge) \cite{globrigid2}. No combinatorial characterization of global rigidity is known for $K>2$. 

\subsubsection{Relevance}
\label{s:rigrelevance}
The combinatorial characterization of (generic) rigidity is important both because of applications --- it helps to know whether a weighted graph is very likely to have finitely or uncountably many realizations --- and because the problem has been open ever since Maxwell's times (mid-1800).

From a practical point of view, however, the computation of the rank of the rigidity matrix, coupled with Asimow and Roth's theorem (Eq.~\eqref{eq:asimowA}-\eqref{eq:asimowB}), appears to settle the question with high probability in polynomial (and practically acceptable) time. This ``reduces'' the problem to one of a mostly theoretical nature. As many theoretical problems, this one remains important also because, in trying to answer this question, researchers keep discovering interesting and practically relevant ideas (such as those about determining global rigidity, see Sect.~\ref{s:globrig}).

Specifically, however, this problem is very important because it has the merit of having shifted the study of rigidity from construction engineering to mathematics. Previously, definitions were rare, ideas informal and therefore ambiguously stated and often confused, and decision procedures almost completely empirical \cite{grunbaum}. Although construction engineers have been building truss structures (such as bridges) for a long time, lack of understanding of concepts such as minimal and redundant rigidity have caused ruptures and disasters throughout history. See the Wikipedia page {\small \url{en.wikipedia.org/wiki/List_of_bridge_failures}}: truss-based bridges that collapsed due to ``poor design'' are likely suspects.

\subsection{Computational complexity of matrix completion problems}
\label{s:complexity}
A common way to define the EDMCP is: given a partially defined matrix, determine whether it can be completed to an sqEDM or not. This problem, introduced in Sect.~\ref{s:matcompl}, is almost like the DGP, except for the seemingly minor detail that $K$ shifts from being part of the input (in the EDGP) to part of the output (in the EDMCP). 

This difference is only apparently minor, as we shall see. While the EDGP is {\bf NP}-hard in the TM computational model, the EDMCP is not known to be hard or tractable for the class {\bf NP}. From a purely computational point of view, the EDGP can only be solved by exponential-time algorithms, while the EDMCP can be solved efficiently, since it can be formulated exactly by the following feasibility-only Semidefinite Programming (SDP) problem \cite{isco16}:
\begin{equation}
\left.\begin{array}{rrcl}
\forall \{i,j\}\in E & B_{ii} + B_{jj} - 2 B_{ij} &=& d_{ij}^2 \\
& B&\succeq& 0.
\end{array}\right\}\label{eq:edmcpsdp}
\end{equation}
The notation $B\succeq 0$ indicates that we require that $B$ is symmetric and positive semidefinite, i.e.~all its eigenvalues must be non-negative.
%% up to page 17 of the corrections
We recall that SDPs can be solved approximately to any desired error tolerance in polynomial time \cite{ipmsdp}. Unfortunately, no method is known yet to solve SDPs {\it exactly} in polynomial time, and membership in the complexity class {\bf P} in the TM computational model is defined by exhibiting a polynomial-time algorithm which can provably tell {\it all} YES instances apart from NO ones.

With respect to the SDP formulation Eq.~\eqref{eq:edmcpsdp}, we remark that, in general, SDP solvers can find solutions of any rank, although a theoretical result of Barvinok \cite{barvinok} shows that if Eq.~\eqref{eq:edmcpsdp} is feasible, then there must exist a solution of rank $K=\lfloor(\sqrt{8|E|+1}-1)/2\rfloor$. Another result of the same author \cite{barvinok2} also proves that, provided there exists a manifold $X$ of solutions of rank $K$, in some sense the solution of Eq.~\eqref{eq:edmcpsdp} ``cannot be too far'' from $X$ in some well defined but asymptotic sense (the latter result has been exploited computationally with some success \cite{isco16}).

This begs the question that is the subject of this section, i.e.~what is the complexity status of the EDMCP? Is it {\bf NP}-hard? Is it in {\bf P}? Is it somewhere between the two classes, provided $\mbox{\bf P}\not=\mbox{\bf NP}$? The rest of the section provides a tutorial towards understanding this question.

\subsubsection{Complexity in the TM model}
\label{s:complbasic}
We first give a short summmary of complexity classes in the TM model. We limit our attention to classes of {\it decision problems}, i.e.~problems to which the answer can be encoded in a single bit. For example, ``given a graph $G$, is it complete?'' is a valid decision problem. Note that the question is parametrized over the input (in this case, the graph $G$): strictly speaking, no solution can be provided until we know what $G$ is. When we replace the symbol $G$ by an actual graph (stored on the TM tape), we obtain an {\it instance} of the decision problem. Thus, decision problems are also seen as infinite sets of instances. The instances having answer bit 1 are known as {\it YES instances}, while the rest are known as {\it NO instances}.

An algorithm $A$ {\it solves} a problem $P$ when it is parametrized by the input $\iota$ of $P$ such that, for each input $\iota$ of $P$, $A(\iota)=1$ if and only if the instance is YES.

The worst-case running time of an algorithm is expressed as a class of functions of the input size $|\iota|=\nu$, i.e.~the amount of bits that $\iota$ needs to be stored on the TM tape. Interesting classes of functions are constants, logarithms, polynomials and exponentials of $\nu$. If the function is $f(\nu)$, the class is indicated as $O(f(\nu))$, and contains all of the functions that are asymptotically upper bounded by $g(\nu) f(\nu) + h(\nu)$ where $g,h$ are themselves asymptotically upper bounded by $f$. A function $g$ is {\it asymptotically upper bounded} by a function $f$ if there is $\nu_0\in\mathbb{N}$ such that for all $\nu>\nu_0$ we have $g(\nu)\le f(\nu)$. For example, if an algorithm takes $36\nu^2+15\nu+3$ CPU cycles to run, then it belongs to the class $O(\nu^2)$. 

For a given problem $P$, it makes sense to ask the asymptotic worst-case running time of the {\it fastest} algorithm that solves $P$ (over all algorithms that solve $P$). For example, answering the question ``given a graph, is it complete?'' admits a trivial algorithm in time $O(|V|^2)$: for each pair of distinct vertices in $V$, check that they are adjacent to an edge. Since there are $|V|(|V|-1)$ pairs, the time is proportional to $|V|^2-|V|$, which belongs to $O(|V|^2)$. We do not know whether this is the best algorithm, but any better algorithm will have running time asymptotically upper bounded by the function $|V|^2$, and so it will belong to the class $O(|V|^2)$.

Since the work of Cobham \cite{cobham} and Edmonds \cite{edmonds} in 1965, we list problems that have a polynomial-time solution algorithm in a class called {\bf P}: this is because any finer granularity would make the algorithm depend on the implementation details of the TM (number of heads, number of tapes and so on), whereas the class of all polynomials is an invariant. The class {\bf P} is known as the class of ``tractable'' problems.

Another interesting class, called {\bf NP}, includes all problems for which YES instances can be {\it proved} to be YES by means of a {\it certificate} that can be verified in polynomial time. In the case of the toy problem ``given a graph, is it complete?'', the certificate is the input itself: as we have seen, the graph can be checked to be complete in polynomial time. This particular feature, i.e.~that the input is the certificate, is shared by all problems in {\bf P}, since the definition of {\bf P} is exactly that there is a polynomial algorithm that, on YES instances, can provide the answer YES in polynomial time. Hence $\mbox{\bf P}\subseteq\mbox{\bf NP}$. The question whether $\mbox{\bf P}=\mbox{\bf NP}$ or not is the most important open question of all computer science and one of the most important in mathematics, and will not be discussed further here (see \cite{johnson1982,aaronson} for more information).

{\bf NP} is an interesting class because it contains problems for which no polynomial time algorithm is currently known, but that are very relevant in practice, such as many packing, covering, partitioning, clustering, scheduling, routing problems, as well as many combinatorial problems with resource constraints: once a solution is given, checking that it solves the problem is generally a question of replacement of symbols by value, and function evaluation, all of which can usually be carried out in polynomial time. On the other hand, finding the solution usually takes exponential time, at least with the algorithms we know so far.

Since it contains so many seemingly hard problems, it makes sense to ask what problems are ``hardest'' in the class {\bf NP}. A qualitative definition of the notion of ``hardest'' is based on the concept of {\it polynomial reductions}. A problem $Q$ can be reduced to another problem $P$ if there is a polynomial time algorithm $\alpha$ mapping YES instances of $Q$ to YES instances of $P$ and NO instances of $Q$ to NO instances of $P$. Now suppose there is a reduction from {\it any} problem $Q$ in {\bf NP} to a given problem $P$. Suppose $P$ were in {\bf P}: then there would exist a polynomial time algorithm $A(\iota)$ to solve each instance $\iota$ of $P$. But then, given an instance $\eta$ of $Q$, the algorithm $A(\alpha(\eta))$ would provide an answer to $\eta$ in polynomial time. In other words, if $P$ were tractable, every problem in {\bf NP} were tractable: even the hardest problems of {\bf NP}. This means that $P$ must be as hard as the hardest problems of {\bf NP}. So we define {\it hardest} for the class {\bf NP} as those problems for which every problem in {\bf NP} can be reduced to them. We call {\bf NP}-hard the class of the hardest problems for {\bf NP}. Note that a problem need not be in {\bf NP} itself in order to be {\bf NP}-hard.

The first such problem was {\sc Satisfiability} (SAT) \cite{cook}: Stephen Cook used a reduction from a polynomial time bounded Turing Machine to a set of boolean constraints of SAT. Ever since, it suffices to reduce {\it from} an {\bf NP}-hard problem $Q$ to a new problem $P$ in order to prove its {\bf NP}-hardness. This can be informally seen as follows: suppose $P$ were {\it not} {\bf NP}-hard, but easier. Then we could combine the solution algorithm for $P$ with the polynomial reduction to yield a proof that $Q$ is not hardest for {\bf NP}, which is a contradiction. The first researcher to spot this reasoning was Richard Karp \cite{karp} in 1972. Now hundreds of problems have been shown to be {\bf NP}-hard \cite{gareyjohnson}. A problem is {\bf NP}-complete if it is both {\bf NP}-hard and belongs to {\bf NP}.

We remark that co{\bf NP} is the class of decision problems that have a polynomial-time verifiable certificate for all NO instances (or, in other words, a polynomial-time verifiable refutation). We remark that {\bf P} is contained in $\mbox{\bf NP}\cap\mbox{co{\bf NP}}$ since the polynomial-time algorithm that decides whether the instance is YES or NO provides, with its own execution trace, a polynomial-time proof (when the instance is YES) as well as a polynomial-time refutation (when the instance is NO). 

\subsubsection{Complexity of {\sc Graph Rigidity}}
\label{s:rigcompl}
It is stated in \cite{recskisurvey2} that the complexity of determining whether a given graph is (generically) infinitesimally rigid in $\mathbb{R}^K$ is in {\bf NP}, since, given a framework $(G,x)$, where the realization $x$ plays the role of certificate, it suffices to compute the rank of the rigidity matrix to establish rigidity according to Eq.~\eqref{eq:asimowA}-\eqref{eq:asimowB}. This assertion is not false, but it should be made clear that it does not refer to the TM computational model. So far, to the best of our knowledge, there is no algorithm for computing the exact rank of a matrix, that is polynomial-time bounded in the TM computational model.

Investigations on the complexity of computing matrix rank do exist, but they are either based on different models of computation, such as real RAM or number of field operations \cite{cheung,burgisser}, or else they place the problem in altogether different complexity classes than {\bf P} or {\bf NP}. Specifically, computing (as well as verifying) the matrix rank over rationals appears to be related to {\it counting} classes. These classes catalog the complexity of the problems of counting the number of solutions of some given decision problems \cite{allender,hoang,mahajansarma}.

On the other hand, it should be clear from Sect.~\ref{s:combcharact1}-\ref{s:combcharact2} that, in the case $K\in\{1,2\}$, the problem of determining generic rigidity of a graph is in {\bf P}. 

\subsubsection{{\bf NP}-hardness of the DGP}
It was shown in \cite{saxe79} that the DGP is {\bf NP}-hard. More specifically, the proof exhibits a reduction from the {\sc Partition} problem (known to be {\bf NP}-complete) to the EDGP with $K=1$ (denoted EDGP${}_1$) restricted to the class of graphs consisting of a single simple cycle.

The {\sc Partition} problem is as follows. Given $n$ positive integers $a_1,\ldots,a_n$ determine whether there exists a subset $I\subseteq \{1,\ldots,n\}$ such that
\begin{equation}
  \sum\limits_{i\in I} a_i = \sum\limits_{i\not\in I} a_i. \label{eq:partition}
\end{equation}
We reduce an instance $a=(a_i\;|\;1\le i\le n)$ of {\sc Partition} to the simple cycle $C=(V,E)$ where $V=\{1,\ldots,n\}$ and $E=\{\{i,i+1\}\;|\;1\le i<n\}\cup\{\{1,n\}\}$. We weigh the edges of $C$ with the integers in $a$, and let $d$ be the edge weight function such that:
\[d_{1,n} = a_1  \qquad \land \qquad \forall 1<i\le n \quad d_{i-1,i}=a_i.\]
Now suppose $a$ is a YES instance of {\sc Partition}: we aim to show that $(C,d)$ is a YES instance of the EDGP${}_1$. Since $a$ is YES, there is an index set $I$ such that Eq.~\eqref{eq:partition} holds. We construct a realization of $(C,d)$ on the real line inductively as follows:
\begin{enumerate}
\item $x_1=0$
\item for all $1<i\le n$, suppose $x_{i-1}$ is known: then if $i\in I$ let $x_i= x_{i-1}+d_{i-1,i}$, otherwise let $x_i=x_{i-1}-d_{i-1,i}$. 
\end{enumerate}
It is obvious by construction that this realization satisfies all of the distances $d_{i-1,i}$ for $1<i\le n$. It remains to be shown that $|x_1-x_n|=d_{1n}$. We assume without loss of generality that $1\in I$ (the argument would be trivially symmetric in assuming $1\not\in I$). We remark that
\[d_{1n}+\sum\limits_{i\in I\smallsetminus\{1\}}\!\!\!\!(x_i-x_{i-1}) = d_{1n}+\sum\limits_{i\in I\smallsetminus\{1\}}\!\!\!\!d_{i-1,i} = \sum\limits_{i\in I} a_i = \sum\limits_{i\not\in I} a_i = \sum\limits_{i\not\in I} d_{i-1,i} = \sum\limits_{i\not\in I} (x_{i-1}-x_i),\]
where the central equality holds by Eq.~\eqref{eq:partition}. Now from the equality between LHS and RHS we have:
\begin{eqnarray*}
  d_{1n}+\sum\limits_{i\in I\smallsetminus\{1\}} (x_i-x_{i-1}) &=& \sum\limits_{i\not\in I} (x_{i-1}-x_i) \\
  \Rightarrow\quad \sum\limits_{1<i\le n} (x_{i-1}-x_i) &=& d_{1n} \\
  \Rightarrow\quad (x_1 - x_2) + (x_2-x_3) + \cdots + (x_{n-1} - x_n) &=& d_{1n} \\
  \Rightarrow\quad |x_1 - x_n| &=& d_{1n},
\end{eqnarray*}
as claimed.

The converse direction aims at showing that if $a$ is a NO instance of {\sc Partition}, then $(C,d)$ is a NO instance of the EDGP${}_1$. We suppose, to aim at a contradiction, that $(C,d)$ is a YES instance of the EDGP${}_1$ instead, and hence that it possesses a realization $x$ on the real line. The argument now traces $x_i$ (as $i$ goes from $1$ to $n$): we let $I$ be the set of indices $i$ such that $x_{i+1}$ is on the right of $x_i$, and claim that $I$ is a solution of {\sc Partition}. More specifically, we let
\[F = \{\{u,v\}\in E\;|\;x_u\le x_v\}\qquad\land\qquad \bar{F}=\{\{u,v\}\in E\;|\;x_u>x_v\},\]
and let $I = \{i<n\;|\;\{i,i+1\}\in F\}$. We have:
\begin{eqnarray*}
  \sum\limits_{\{u,v\}\in F} (x_u-x_v) &=& \sum\limits_{\{u,v\}\in\bar{F}} (x_v-x_u)\\
\Rightarrow\quad  \sum\limits_{\{u,v\}\in F} |x_u-x_v| &=& \sum\limits_{\{u,v\}\in\bar{F}} |x_v-x_u|\\
\Rightarrow\quad  \sum\limits_{\{u,v\}\in F} d_{uv} &=& \sum\limits_{\{u,v\}\in\bar{F}} d_{uv}\\
\Rightarrow\quad  \sum\limits_{i\in I} a_i &=& \sum\limits_{i\not\in I} a_i,
\end{eqnarray*}
against the assumption that $a$ is a NO instance of {\sc Partition}. 

The above argument shows that the subclass EDGP${}_1$ of the DGP, restricted to simple cycle graphs, is {\bf NP}-hard. Since it is contained in the whole DGP class, then the DGP itself must be {\bf NP}-hard, for if it were not, then it would be possible to solve even the restricted subclass in an easier way.

Hardness for any complexity class is defined using appropriate reductions. 

\subsubsection{Is the DGP in {\bf NP}?}
In general, the DGP is not known to be in {\bf NP}. The reason is that, in the TM model, inputs are always rational numbers, which means that Eq.~\eqref{eq:dgp} might have irrational solutions: in other words, the realizations of the given graph could involve irrational (though algebraic) numbers in the components. Although there are a few finitistic encodings of algebraic numbers, none of them is a good candidate to verify that the realization satisfies Eq.~\eqref{eq:dgp} exactly \cite{dgpinnp}. The EDMCP is not in {\bf NP} for much the same reasons. Both, however, are in {\bf NP} with respect to the real RAM computational model \cite{blum}. 

We remark that the EDGP${}_1$ is in {\bf NP}: any irrational realization can arbitrarily be translated to be aligned with a rational point for any vertex (say vertex $1$). Since the distances are all rational (as they are part of the input), all of the other vertices will have rational points too. These rational points can be verified to satisfy the square of Eq.~\eqref{eq:dgp} exactly in polynomial time. Since the EDGP${}_1$ is both {\bf NP}-hard and in {\bf NP}, it is {\bf NP}-complete.

It is shown in \cite{oneinfnorm} that the DGP in $\ell_1$ and $\ell_\infty$ norms belongs to {\bf NP} independently of its dimension. 

\subsubsection{EDMs and PSD matrices}
\label{s:edmpsd}
As mentioned in Sect.~\ref{s:history} above, EDMs and Positive Semidefinite (PSD) matrices are strongly related. Let $D$ be a sqEDM, and assume without loss of generality that $D$ is yielded by an $n\times K$ realization matrix that is centered, i.e.~$\sum\limits_{v\le n} x_v = \mathbf{0}_K$ (the all-zero $K$-vector). Then from the identity:
\begin{equation}
 \forall u,v\le n\quad \|x_u-x_v\|_2^2 = \|x_u\|_2^2 + \|x_v\|_2^2 - 2 x_u\cdot x_v, \label{eq:eucldist}
\end{equation}
we obtain the matrix equation:
\begin{equation}
  D = r\,\mathbf{1}^\top + \mathbf{1}\,r^\top -2 B, \label{eq:mateucldist}
\end{equation}
where $B=x\,x^\top$ is the {\it Gram matrix} of $x$, $r$ is the column $n$-vector consisting of the diagonal element of $x\,x^\top$, and $\mathbf{1}$ is the all-one $n$-vector. 

We now consider the {\it centering matrix} $J=I - \frac{1}{n}\mathbf{1}\,\mathbf{1}^\top$, where $I$ is the $n\times n$ identity matrix: when applied to an $n\times K$ matrix $y$, $Jy$ is translation congruent to $y$ such that $\sum_{i\le n} (Jy)_i = \mathbf{0}_K$. We remark that $J$ is symmetric, so $J^\top=J$. From Eq.~\eqref{eq:mateucldist} we obtain:
\begin{eqnarray*}
  -\frac{1}{2}\, J\, D\, J &=& -\frac{1}{2} (J\,r\,\mathbf{1}^\top\, J + J\,\mathbf{1}\,r^\top\,J) + J\,B\,J \\
  &=& -\frac{1}{2}(J\,r\,\mathbf{0}_n^\top + \mathbf{0}_n\,r^\top\,J) + B = B,
\end{eqnarray*}
whence
\begin{equation}
  B = -\frac{1}{2}\,J\,D\,J. \label{eq:BDsq}
\end{equation}
Note that $J\,\mathbf{1}=\mathbf{0}_n$ since centering the all-one vector trivially yields the all-zero vector, and that $J\,B\,J=J\,x\,x^\top\,J=x\,x^\top=B$ since $x$ was assumed centered at zero.

It is easy to see that a matrix is Gram if and only if it is PSD. If $B$ is a Gram matrix, then there is an $n\times K$ realization matrix $x$ such that $B=x\,x^\top$. Let $\bar{x}$ be the $n\times n$ matrix obtained by padding $x$ on the right with zero columns, then $B=\bar{x}\,\bar{x}^\top$ is a valid (real) factorization of $B$, which means that all eigenvalues of $B$ are non-negative. Conversely, suppose $B$ is PSD; the eigendecomposition of $B$ is $P^\top\,\Lambda\,P$, where $\Lambda$ is a diagonal matrix with the eigenvalues along the diagonal. Since $B$ is PSD, $\Lambda\ge 0$, which means that $\sqrt{\Lambda}$ is a real diagonal matrix. Then by setting $x=P^\top\,\sqrt{\Lambda}$ we have $B=x\,x^\top$ which proves that $B$ is a Gram matrix.

These two results prove that $D$ is a sqEDM if and only if $B$ is a PSD matrix. By this equivalence, it follows that the EDMCP is in {\bf P} if and only if the {\sc PSD Completion Problem} (PSDCP) is in {\bf P}. Given the wealth of applications that can be modelled using SDP, establishing whether the PSDCP is in {\bf P} makes the same question for the EDMCP even more important.

In \cite{laurent00}, it is observed that the PSDCP belongs to {\bf NP} if and only if it belongs to co{\bf NP}, so the PSDCP cannot be either {\bf NP}-complete or co{\bf NP}-complete unless $\mbox{\bf NP}=\mbox{co{\bf NP}}$.

\subsubsection{Relevance}
\label{s:complrelevance}
Why is it important to determine the complexity of the EDMCP (or, equivalently, of the PSDCP)? Similarly to what was said in Sect.~\ref{s:rigrelevance}, from a practical point of view, we can solve EDMCPs and PSDCPs in polynomial time to any desired accuracy using SDP solver technology. On the other hand, this is only a partial answer.

Theoretically speaking, the issue as to whether $\mbox{\bf P}=\mbox{\bf NP}$, $\mbox{\bf P}\not=\mbox{\bf NP}$ or even whether the whole question might actually be independent of the ZFC axioms \cite{aaronson} dominates the field of theoretical computer science and much of discrete mathematics. Identifying a problem that is in $\mbox{\bf NP}\smallsetminus\mbox{\bf P}$ would obviously settle the question in the negative. Every such candidate is interesting until eliminated. Even though the EDMCP is not the perfect candidate on account of the uncertain status of its membership in {\bf NP}, it is still one of the few known\footnote{Another, and possibly better, candidate is the {\sc Graph Isomorphism} problem \cite{babai,babai2}.} practically relevant problems with this status.

Another, and more practical reason is the bottleneck represented by current SDP solver technology: although the algorithm is polynomial-time, and notwithstanding the fact that existing implementations are of high quality \cite{mosek7}, solving SDPs with significantly more than thousands of variables and constraints is still a hurdle. It is hoped that research in complexity will drive a search for better SDP solution algorithms.

\subsection{Number of solutions}
\label{s:numsol}
The issue of finding, or estimating the number of solutions of a DGP instance prior to solving it was brought forward as a ``major challenge in DG'' by the scientific committee of the DGTA16 workshop \cite{dgta16proc} when putting together a (successful) NSF proposal to support the workshop.

Although counting solutions of DGP instances is related to establishing rigidity or flexibility of graphs and frameworks (see Sect.~\ref{s:rigidity}), the former is a finer-grained alternative to the latter. Rigidity results are qualitative insofar as they focus on three categories: unique solutions, a finite number of solutions, and infinitely many solutions. For some applications, this is not enough, and the exact or approximate number of solutions is required, or at least helps in the solution process. In this sense, counting solutions deserves the status of {\it open research area} independently of the studies on rigidity.

\subsubsection{Either finite or uncountable}
\label{s:finuncount}
The first question that might arise is: can any DGP instance have an infinite, but countable number of solutions? The answer to this question is negative, and comes from both both topology and algebraic geometry.

The ``topological'' proof rests on an observation made by John Milnor in 1964 \cite{milnor64}. If $\mathcal{V}\subseteq\mathbb{R}^m$ is the variety defined by the system of $p$ polynomial equations
\begin{equation}
  \forall i\le p \quad f_i(x_1,\ldots,x_m) = 0
  \label{eq:variety}
\end{equation}
each of which has degree bounded above by the integer $k\ge 1$, then the sum of the Betti numbers of $\mathcal{V}$ is bounded above by $k(2k-1)^{m-1}$.

A straight application of this result to the squared version of Eq.~\eqref{eq:dgp} encoding the DGP leads to $k=2$, $m=Kn$, $\mathcal{V}$ being the set of realizations satisfying the given DGP instance, and the sum of the Betti numbers being bounded above by $2(3^{Kn-1})$, which is $O(3^{Kn})$. We now recall that the Betti numbers count, among other things, the connected components of a variety. This immediately yields that that there are finitely many connected components. It is well known from basic topology that connected components can either be uncountable sets or consist of an isolated point.

A second method of proof consists in invoking the cylindrical algebraic decomposition results in real algebraic geometry \cite{benedetti,pollack}, which stem from Tarski's quantifier elimination theory \cite{tarski-reals}, to show that $\mathcal{V}$ consists of a finite number of connected components with a certain (cylindrical) structure. The result follows.

While the algebraic geometry result is quantitative, meaning that it provides an explicit description of the geometrical properties of each connected component, the topological method simply proves the finiteness of the number of connected components, which is really all that is required.

\subsubsection{Loop flexes}
\label{s:loopflex}
\begin{figure}[!ht]
  \begin{center}
    \includegraphics[width=7cm]{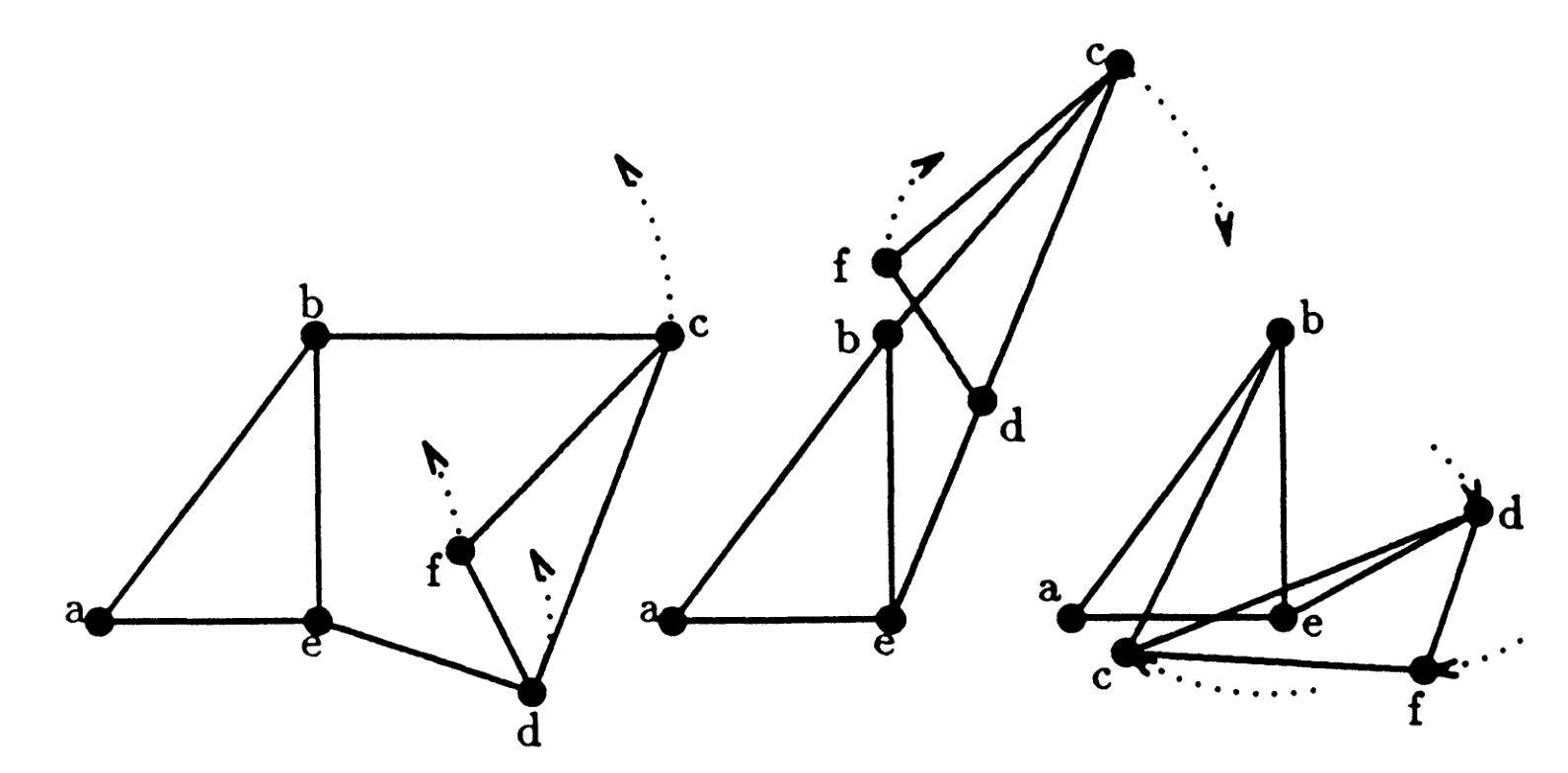}
  \end{center}
  \caption{Each flex is diffeomorphic to a circle.}
  \label{f:hen92}
\end{figure}

An interesting observation about the ``shape'' of flexes in flexible graphs was made by Bruce Hendrickson in 1992 \cite[Thm.~5.8]{Hen92}: if a graph is connected, (generically) flexible in $\mathbb{R}^K$ and has more than $K+1$ vertices, then for almost all edge weight functions the realization set contains a submanifold that is diffeomorphic to a circle. The situation is sketched graphically in \cite[Fig.~7]{Hen92}, reported in Fig.~\ref{f:hen92}.

A {\it manifold} is a topological set which is locally homeomorphic to a Euclidean space; a {\it homeomorphism} is a continuous function between topological spaces which also has a continuous inverse function. A {\it diffeomorphism} is a smooth invertible function that maps a differentiable manifold to another, such that the inverse is also smooth. A function is {\it smooth} if it has continuous derivatives of any order everywhere. And, lastly, a manifold is {\it differentiable} if the composition with a local homeomorphism with the inverse of any other local homeomorphism is a differentiable map. All of this differential topology definitions formalize the concept that, although none of the graph vertices might really move in a circle, the flex itself contains a closed loop that is topologically equivalent to a circle.

The proof proceeds by adding to the flexible graph $G$ as many edges as are required to leave just one degree of freedom to the flex. It is then relatively easy to show that the flex is compact and that it is a manifold. A well known result\footnote{See e.g.~\url{en.wikipedia.org/wiki/Classification_of_manifolds}.} in differential topology shows that compact one-dimensional manifolds are diffeomorphic to circles.

This result was used by Hendrickson to prove that redundant rigidity is a necessary condition to solution uniqueness: if a graph is rigid but not redundantly so, then there is an edge $\{i,j\}$ such that its removal yields a flex diffeomorphic to a circle. In almost all cases, there will be two points on this circle, corresponding to two incongruent realizations, that are compatible with the given edge distance $d_{ij}$ (in Fig.~\ref{f:hen92}, the missing edge $\{a,f\}$ has the same length $\|x_a-x_f\|_2$ in both the left and the right picture).

\subsubsection{Solution sets of protein backbones}
Proteins are organic molecules at the basis of life. They interact with other proteins and with cells and living tissues by chemical reaction of the atoms in their outside shell. For these reactions to happen, the protein must physically be able to ``dock'' to the prescribed sites, which involve it having a certain geometrical shape. This way, proteins activate and inhibit living processes in all animals. Proteins consist of a backbone together with some side chains \cite{sidechains}, the building block of which are a set of around twenty small atomic conglomerates called {\it amino acids} \cite{wuthrich_83}. The problem of finding the shape of the protein can be decomposed in finding the shape of the backbone and then placing the side chains correctly \cite{santana07,santana08}.

The backbone itself has an interesting graph structure, in that it provides an atomic order with a very convenient geometric property: we can estimate rather precisely the distances between each atom having rank greater than two in the order and its two immediate predecessors. We know the distances $d_{i-1,i}$ for each $i>1$ because they are covalent bonds; and, since we also know all of the covalent angles, we can also compute all of the missing sides from the consecutive triangles, which yields the distances $d_{i-2,i}$. This yields a graph consisting of $i-2$ consecutive triangles adjacent by sides. Since most protein graphs need realizations in $\mathbb{R}^3$, the sides by which the triangles are attached provide some rotating hinges, which in turn means that we have uncountably many solutions.

\begin{figure}[!ht]
  \begin{center}
    \includegraphics[width=10cm]{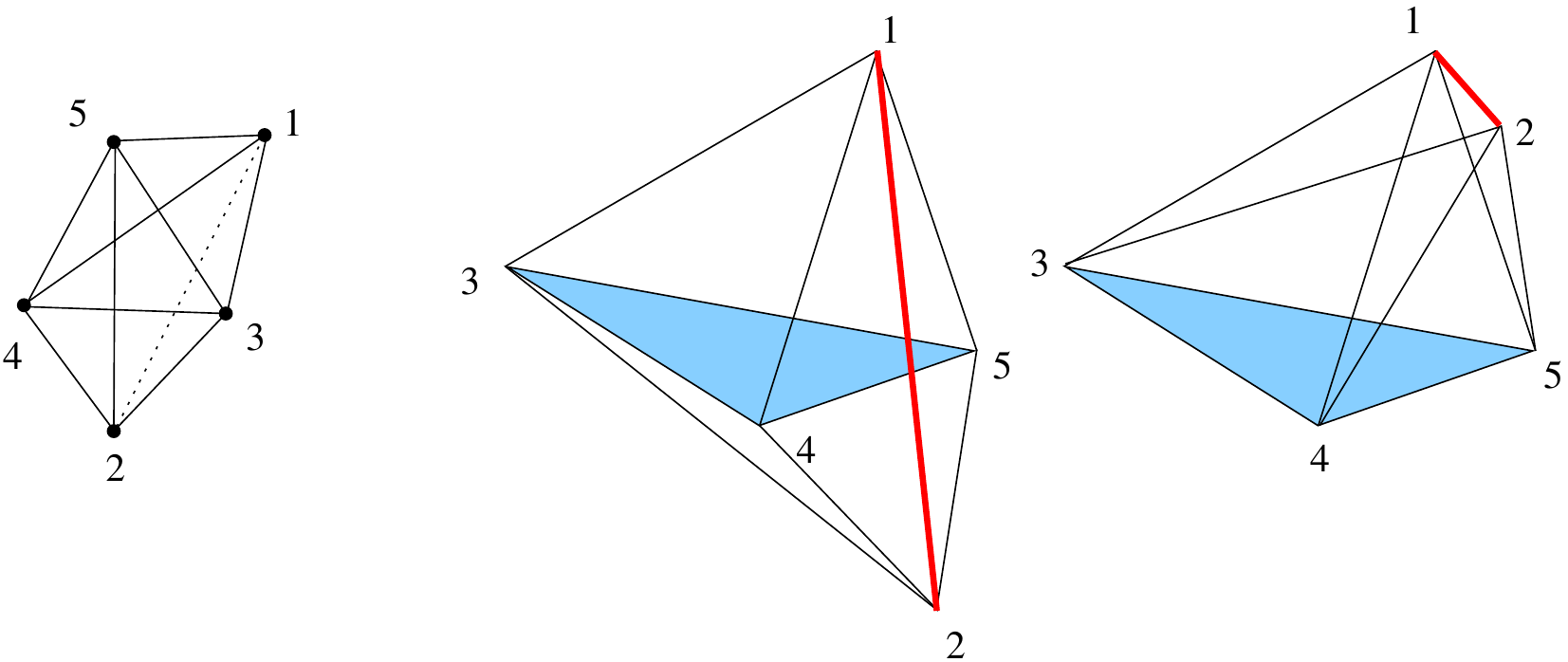}
\end{center}
\caption{\footnotesize Two realizations of a quasiclique in $\mathbb{R}^3$
  (the missing edge is dotted on the left, and bold on the center and
  on the right).}
\label{f:realquasiclique}
\end{figure}

This is where W\"uthrich's Nobel prize winning NMR based techniques come in: we can estimate most distances smaller than about 6{\AA}. Note that generally, the distances $d_{i-3,i}$ are smaller than this threshold. After adding these distances to our protein graph, we have a structure consisting of a sequence of consecutive tetrahedra adjacent by faces, which is clearly rigid in 3D, though not globally rigid, as shown in Fig.~\ref{f:realquasiclique}. 
From the point of view of the underlying graph, protein backbones are defined by containing a subgraph consisting of a sequence of consecutive $4$-cliques adjacent by $3$-cliques. When such a sequence consists of two $4$-cliques it is also known as a {\it $5$-quasiclique}, since it is a graph on 5 vertices with all edges but one.

\begin{figure}[!ht]
\begin{center}
\includegraphics[width=5cm]{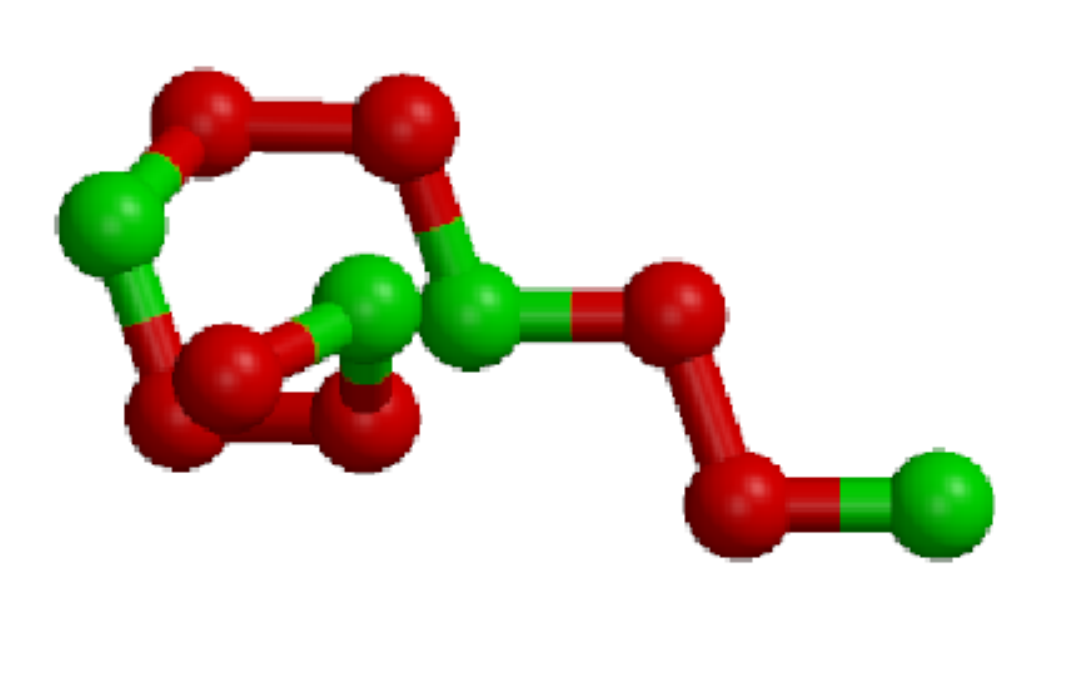}
\end{center}
\caption{\footnotesize A 3D realization of an artificial protein backbone \cite{Lav05}.}
\label{f:backbone}
\end{figure}

Recall that NMR estimates distances up to 6{\AA}: while these contain those between atoms $i$ and $i-3$, they may also contain other distances if the protein backbone folds back in space close to itself (see Fig.~\ref{f:backbone}).
So the edge set $E$ of a typical protein backbone graph $G=(V,E)$ also has these other distances. We call these {\it pruning distances} (and the underlying edges {\it pruning edges} $E_P$), while the distances forming the clique sequence are called {\it discretization distances} (and the underlying edges {\it discretization edges} $E_D$).

The subset of the MDGP containing these instances is called the {\sc Discretizable MDGP} (DMDGP) \cite{dmdgp}. This problem was shown to be {\bf NP}-hard in \cite{bppolybook}. It can be solved using an algorithm called Branch-and-Prune (BP) \cite{lln5,mdjeep}, which works inductively as follows: suppose atom $i-1\ge 3$ has already been realized as $x_{i-1}\in\mathbb{R}^3$. Then, by trilateration \cite{eren04}, with probability 1 there will be at most two positions for $x_i$.

\begin{figure}[!ht]
  \begin{center}
    \includegraphics[width=4cm]{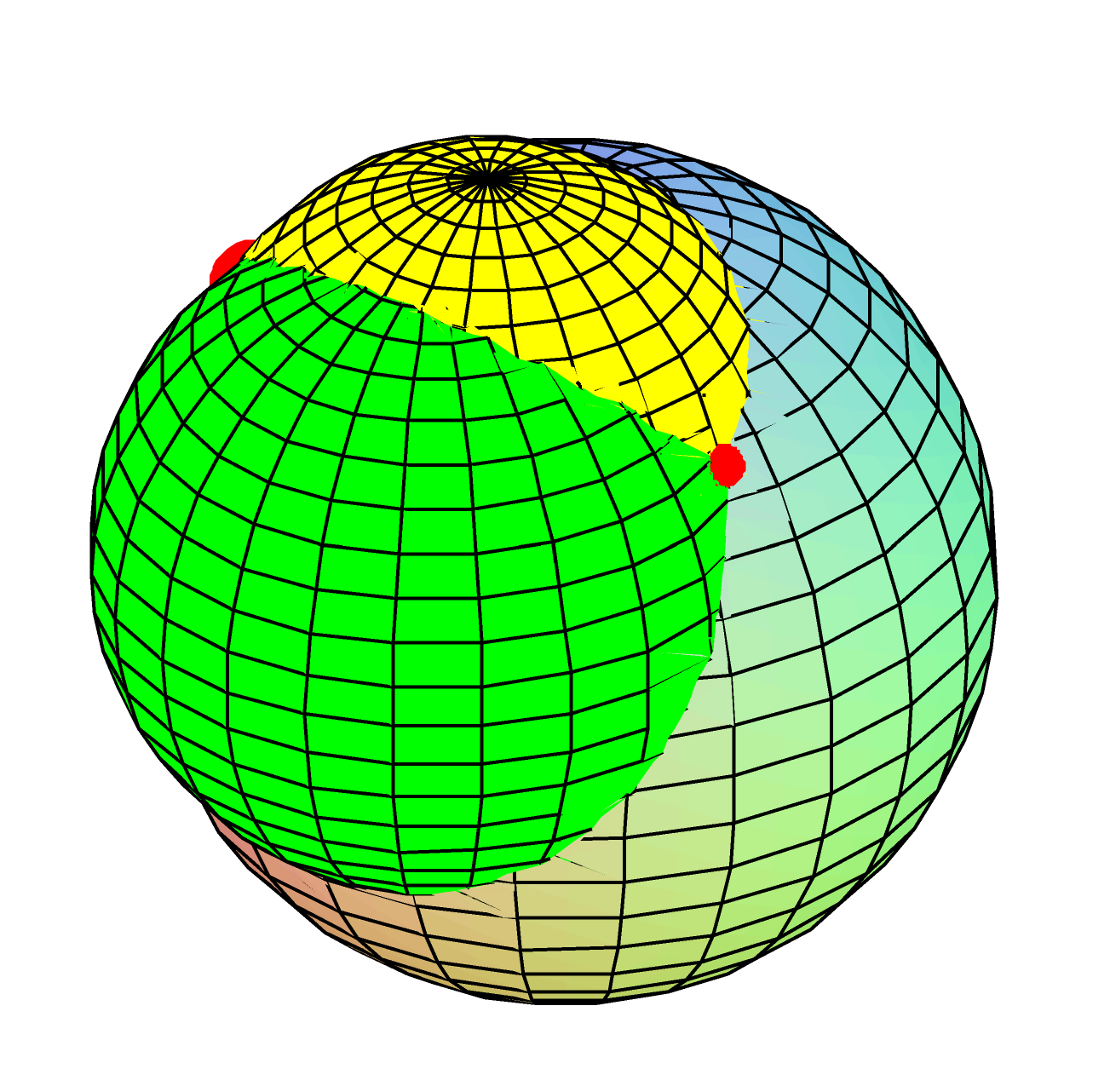}
  \end{center}
  \caption{The intersection of three spheres in $\mathbb{R}^3$ contains at most two points.}
  \label{f:3spheres}
\end{figure}

This can be seen e.g.~in Fig.~\ref{f:realquasiclique} using the order $(1,3,4,5,2)$: once vertex $5$ is realized, there are two positions for vertex $2$, one close to vertex $1$ (Fig.~\ref{f:realquasiclique}, right) and one further away (Fig.~\ref{f:realquasiclique} center). The probability zero cases, ignored by the algorithm, are those where the distances are exactly right for vertex $2$ to be realized coplanar with $x_3,x_4,x_5$, in which case there is only one position for vertex $2$. Another intuitive way of seeing this fact is that vertex $i$ is at the intersection of the three spheres $S(x_{i-3},d_{i-3,i})$, $S(x_{i-2},d_{i-2,i})$, $S(x_{i-1},d_{i-1,i})$. Such an intersection contains at most two points \cite{coope,dvop} (see Fig.~\ref{f:3spheres}).

Once two positions have been found for vertex $i$, say $x_i^+,x_i^-$, the BP algorithm checks whether either, both or neither of these are compatible with any pruning distance which might be adjacent to vertex $i$: the incompatible positions are pruned (ignored). After this step, there may remain zero, one or two compatible positions. If there are none, this particular branch of the recursion is pruned, and the search backtracks to previous recursion levels. If there is only one position, it is accepted, and the search is called recursively. If there are two positions, the search is branched, and each branch yields a recursive call. When $x_n$ is placed (where $n=|V|$), a realization $x$ has been found: this recursion path ends and $x$ is stored in a solution set $X$. The search is then backtracked, until no more recursion calls are possible.

The BP algorithm yields a search tree of height $n$ with several interesting properties.
\begin{enumerate}[(a)]
\item A necessary branching occurs at $x_4$ \cite{dmdgp} since no pruning edges involving predecessors of $4$ can be adjacent to it.\label{enum:4lev}
\item If the pruning edges contain all of the $\{i-4,i\}$ pairs (for $i>4$) then the only branching that occurs is the ``$4$-th level branching'' referred to in (\ref{enum:4lev}) above, and the BP runs in worst-case polynomial time.
\item If the pruning edges contain $\{1,n\}$ then the instance has only two realizations which are in fact reflections of each other through the symmetry induced by the $4$-th level branching. Notwithstanding, BP does not necessarily run in polynomial time since extensive branching might occur, only to be pruned before or at the last ($n$-th) level.
\item The worst-case running time of BP, in general, is exponential in $n$. In practice, however, it is very fast for considerably large instances, and very precise.
\item The BP can be stopped (heuristically) based on the number of solutions found, to make it faster. But it is still {\bf NP}-hard to find even one solution, so it could behave exponentially nonetheless.
\item There is nothing special about three-dimensional spaces: the BP also works in arbitrary Euclidean spaces $\mathbb{R}^K$ for any $K\ge 1$.
\item When run to completion, in general the BP finds all solutions modulo rotations and translations; if one of the two sides of the $4$-th level branching is pruned, then $X$ will contain all incongruent realizations for the given instance.
\item The pruning edges induce a very elegant pattern of partial reflections over the backbone. Let $T=\bigcup X$ be the superposition in space of all of the realizations in $X$. It is shown in \cite{powerof2,bppolybook} that the set $T$ is invariant to a certain set of {\it partial reflection} operators $g_i$, acting on realizations, that fix all of the vertex positions up to $x_{i-1}$, and then reflect $x_i,\ldots,x_n$ with respect to the affine subspace spanned by $x_{i-1},\ldots,x_{i-K}$. Moreover, the partial reflection group can be computed {\it a priori} as
  \begin{equation}
    G_P = \langle g_i\;|\;i>K\land \not\exists \{u,v\}\in E_P\;(u+K<i\le v)\rangle. \label{eq:pruninggroup}
  \end{equation}
\item This symmetry structure makes it possible to compute a single realization $x$, and then generate $X$ as the orbit $G_Px$ of $x$ \cite{symmBPjbcb}.
\item More relevant to this open research area, this symmetry structure also allows the aprioristic computation of the exact number of partial realizations active each of the $n$ levels of the BP search tree \cite{liberti-gsi13}. In particular, this allows the control of the width (i.e., the breadth) of the tree in function of the distribution of the pruning edges \cite{bppolybook}.
\item The estimation of the number of partial solutions at each level of the tree has yielded a crucial empirical observation: if protein backbones fold and the folds change direction only a logarithmic number of times in $n$, then the width of the tree is bounded above by a polynomial in $n$, which implies that the BP is a polynomial time algorithm on such instances. 
\end{enumerate}

The BP actually works on larger classes than the DMDGP: the largest class of instances that the BP can solve is called {\sc Discretizable DGP} (DDGP) \cite{ddgp} and consists of instances having a vertex order such that any vertex $i$ is adjacent to at least three adjacent (not necessarily immediate) predecessors. In this setting, however, the results about partial reflection symmetries listed above are invalid. The BP was also extended to solve instances of the {\it interval} DDGP (iDDGP) problem, i.e.~some of the distances in a DDGP instance are represented as intervals \cite{bpinterval,bipbip}. In this setting, the BP ceases to be an exact and exhaustive algorithm, but it can still find a number of incongruent solutions to the instance.

\subsubsection{Clifford algebra}
\label{s:cliffordalgebra}
A particular type of geometric algebra, called {\it Clifford algebra}, is used to represent compactly and perform computations with various geometrical shapes in Euclidean spaces. It was recently used to provide a compact representation of the solution set of DDGP instances \cite{cliffordalgebra}.

\begin{figure}[!ht]
  \begin{center}
    \includegraphics[width=7cm]{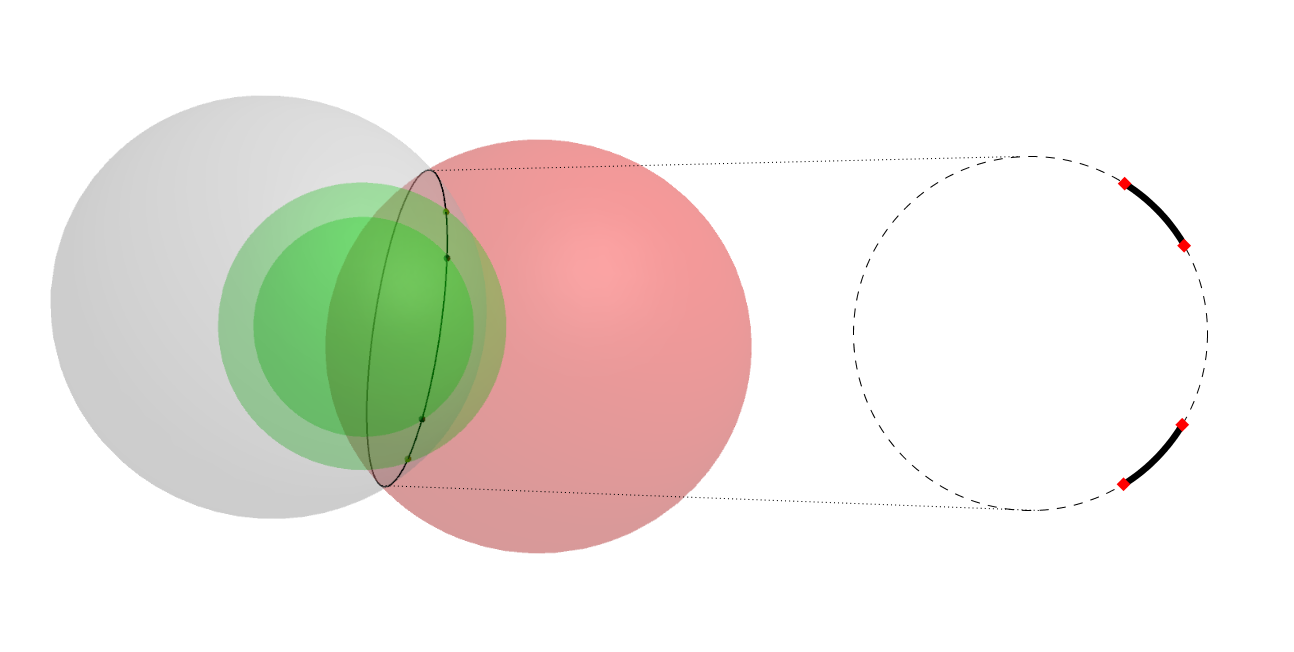}
  \end{center}
  \caption{Intersection of two spheres with a spherical shell: at most two circular arcs.}
  \label{f:sphshell}
\end{figure}

Carlile Lavor's presentation at the DGTA16 workshop \cite{dgta16proc} focused on an interesting extension of this representation to those iDDGP${}_3$ instance where, for each $v>3$, at most one in three distances $d_{uv}$ where $u$ is an adjacent predecessor of $v$ is represented as an interval. In such cases, the intersection of three spheres shown in Fig.~\ref{f:3spheres} becomes the intersection of two spheres and a spherical shell (see Fig.~\ref{f:sphshell}), which turns out to consist of two circular arcs with probability 1.
Currently, efforts are under way to use these representations in order to do branching on circular arcs \cite{clifford-dga13} rather than discretizations thereof \cite{bpinterval}. 

\subsubsection{Phase transitions between flexibility and rigidity}
\label{s:phasetrans}
Consider a random process by which an initially empty graph is added edges with some probability $p$:
\begin{enumerate}
\item sample $u,v$ randomly from $V$
\item if $\{u,v\}\not\in E$, add $\{u,v\}$ to $E$ with probability $p$
\item repeat.
\end{enumerate}
At the outset, the graph is likely to consist of isolated edges (i.e.~paths consisting of one edge), or pairs of adjacent edges (i.e.~paths consisting of two edges), and is therefore flexible in $\mathbb{R}^2$. As more and more edges get added, some rigid components appear that can move flexibly with respect to one another, and finally the whole graph becomes rigid. What is the value of $\eta=\frac{2|E|}{|V|(|V|-1)}$ which marks the appearance of a single rigid component? Similar graph generation processes have been analysed with respect to the edge generation probability for the appearance of giant connected components in random graphs \cite{erdos-renyi,bollobasrndgph}.

In the case of rigidity, Duxbury and Thorpe independently (but in two papers published consecutively in Phys.~Rev.~Lett., in which one cites the other) proposed {\it percolation analysis} in 1995 \cite{duxbury_perc,thorpe_perc} (also see \cite[\S 4]{dg-4or}). Essentially, they randomly add edges (with probability $p$) to an empty graph, and verify rigidity of clusters after each edge addition using the so-called ``pebble game'' (see \url{linkage.cs.umass.edu/pg/pg.html}). This is a graph labelling algorithm which identifies planar isostatisticity (by checking Laman's conditions (a)-(b) in Sect.~\ref{s:combcharact2}) while flagging redundant rigidity, or determines flexibility. These simulations link the emergence of a giant rigid component to the parameters of the random edge addition process, such as e.g.~the edge creation probability $p$. 

In his lectures on rigidity given at the Institut Henri Poincar\'e in 2005 in Paris, Bob Connelly observes that on the type of triangular plane tessellation graphs used by Duxbury and Thorpe, the percolation simulations agree with the theoretical observation that the whole graph needs $p\ge\frac{2}{3}$ in order to be rigid, which is required to achieve $|E|\ge 2|V|-3$ on average. 

An interesting open question about phase transitions is motivated by the following observation: while the DMDGP is {\bf NP}-hard, the DGP${}_K$ (where $K=3$) on protein backbone graphs where the edge $\{i-K-1,i\}$ has been added to the graph for all $i>K+1$ can be solved efficiently using $K$-lateration (moreover, further increasing the number of edges helps, as the DGP can be solved efficiently on complete graphs too \cite{dongwu}). What is the critical value of the parameter $\eta$ (defined above) that determines the phase transition from {\bf NP}-hardness to tractability?

%% correcting Carlile's comments (page 30)

\subsubsection{Relevance}
\label{s:numsolrelevance}
This open area is mostly motivated by applications. There are applications, such as clock synchronization, wireless networks, autonomous underwater vehicles and others, where one is interested in DGP instances with exactly one solution (modulo congruences). There are other applications, such as those to molecular conformation (be it proteins or nanostructures) where one is interested in all (finitely many) chiral isomers of the molecule in question. And there are areas such as robotics, where one is interested in the whole solution manifolds, including flexes.

This is not all. Although we consider molecular graphs to be rigid, which helps with computation, it is well known that atoms vibrate in molecules according to many factors, including the temperature. This means that the molecules undergo internal movement. Although these are not all necessarily of the flex type (meaning that some of them may not preserve all pairwise distances), the strongest chemical bonds may well be (almost preserved). So there is an analysis of flexibility involved \cite{sljoka}. Moreover, there are very few applications where the distances are known precisely. Mostly, distance errors are modelled as intervals, as in the {\it i}DGP (see Sect.~\ref{s:idgp}). This necessarily introduces some flexibility in the frameworks. 

In most of the situations where a DGP instance of interest has more than one solution, it helps to have an {\it a priori} estimation of the number of solutions. In the presence of flexes, it also helps to have an idea of the type of flex involved. 

\subsection{The unassigned DGP}
\label{s:udgp}
The uDGP was briefly introduced in Sect.~\ref{s:unassigned}, but we formally state it here for clarity.
\begin{quote}
  {\sc unassigned DGP} (uDGP). Given a positive integer $K>0$, a graph $G=(V,E)$ and a sequence $D$ of $m=|E|$ positive real values $(d_i\;|\;i\le m)$, determine whether there exists an assignment function $\alpha:\{1,\ldots,m\}\to E$ and a realization $x:V\to\mathbb{R}^K$ such that:
  \begin{equation*}
    \forall \{u,v\}\in E\; \exists i\le m \quad \|x_u-x_v\| = d_{\alpha(i)}.\label{eq:udgp}
  \end{equation*}
\end{quote}
In general, this is a problem schema valid for any norm; usually, the uDGP is of interest in the Euclidean norm.

The uDGP was previously studied only in the context $K=1$ (see e.g.~\cite{skiena,dakic}). It was formally introduced in the context of nanostructure determination (see Sect.~\ref{s:nano}) in \cite{tribond}, as the optimization problem:
\begin{equation}
\min\limits_{\alpha,x} \sum\limits_{\{u,v\}\in E} f(\|x_u-x_v\|_2 - d_{\alpha^{-1}(u,v)}),\label{eq:udgpcost}
\end{equation}
for any strictly convex univariate function $f$ achieving its global minimum at 0.

\subsubsection{Determining nanostructures from spectra}
\label{s:nano}
In his talk at the DGTA16 workshop \cite{dgta16proc}, Simon Billinge jokingly explained that, while the crystal structure determination is largely a solved problem, on account of one giving it to one's grad student so she can push the ``start'' button on the X-ray machine, the {\it nanostructure determination problem} is far from being in the same class.

\begin{figure}[!ht]
  \begin{center}
    \includegraphics[width=7cm]{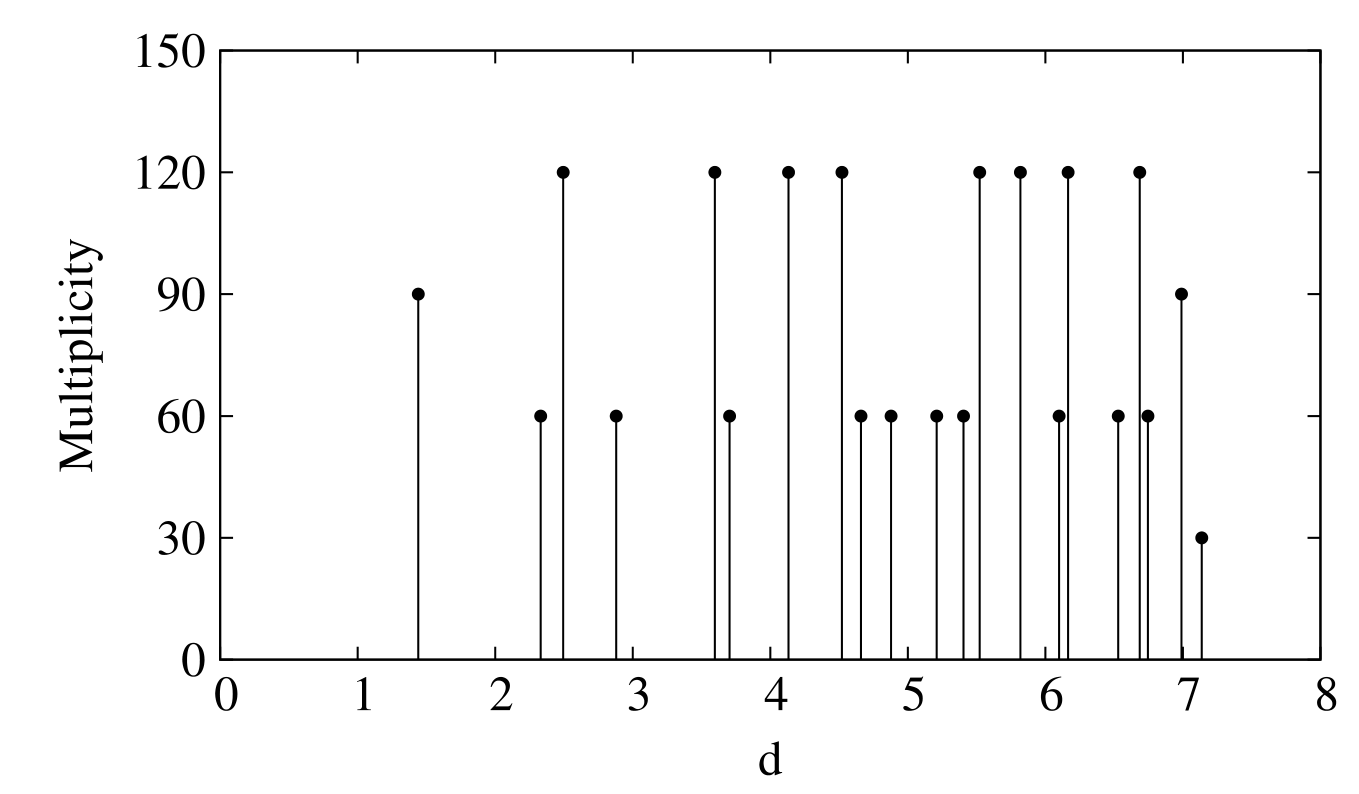}
  \end{center}
  \caption{Toy example of a PDF. PDFs arising from experimental data are noisy and look like continuous wiggly curves with peaks corresponding to observed distances.}
  \label{s:pdf}
\end{figure}

In crystals, the translational symmetry implies that X-ray experiments yield a periodic response signal, which can be decomposed using Fourier analysis. For nanostructures, one can get a response signal only from a large set of similar nanoparticles with unknown orientation. The output of such experiments is a Pair Distribution Function (PDF) $g:\mathbb{R}_+\to[0,1]$, which is a function mapping a distance value $d$ to the frequency with which $d$ occurs in the nanostructure (see Fig.~\ref{s:pdf}).

The specificity of the uDGP input coming from the application to nanostructure determination is that it allows the estimation of {\it all} of the inter-atomic distances \cite{tribond,dg-4or} (so $G$ is complete), and the output realization is required to be in $\mathbb{R}^3$. We recall that, although realizing complete graphs is a tractable problem (either using trilateration \cite{dongwu} or matrix factoring via the connection between EDM and PSD matrices, see Sect.~\ref{s:edmpsd}), the uDGP has the added difficulty of finding the assignment $\alpha$ at the same time as the realization, so it is not clear at all whether the uDGP on complete graphs might be tractable in full generality. 

An algorithm for solving the uDGP on complete graphs, called {\it tribond}, was proposed in \cite{tribond} in full generality for any $K$, and shown in Alg.~\ref{a:tribond}.
\begin{algorithm}
\begin{algorithmic}
  \STATE {\bf Input}: an integer $K>0$, a graph $G=(V,E)$, a sequence of $m=|E|$ values $D=(d_1,\ldots,d_m)$.
  \FOR{each subsequence $\bar{D}$ of $K+2$ distances in $D$}
  \IF{$\bar{D}$ can be realized by a partial realization $x$}
  \FOR{$K+2< i\le n$}
  \FOR{each subsequence $S$ of $K+1$ distances in $D\smallsetminus\bar{D}$}
  \IF{$\not\exists\; x_i\in\mathbb{R}^K$ such that $(x,x_i)$ is a realization for $(\bar{D},S)$}\label{a:tribondif}
  \STATE {\bf break}
  \ELSE
  \STATE $x\leftarrow(\bar{x},x_i)$, $\bar{D}\leftarrow(\bar{D},S)$
  \IF{$i=n$}
  \RETURN{$x$}
  \ENDIF
  \ENDIF
  \ENDFOR
  \ENDFOR
  \ENDIF
  \ENDFOR
  \RETURN{{\tt infeasible}}
\end{algorithmic}
\caption{The tribond algorithm.}
\label{a:tribond}
\end{algorithm}
It is claimed in \cite{tribond} that if the initial guess for $\bar{D}$ is good, and the subsequent choices of vertex subsets $S$ (see Alg.~\ref{a:tribond}) can always be extended to a full solution $x$, then this algorithm runs in polynomial time, since it does not incur the computational cost of looping over all subsets of cardinality $K+2$ and $K+1$. However, it is easy to notice that, if $K$ is fixed, the tribond algorithm always takes polynomial running time. This is because $K+2\choose n$ and $K+1\choose n-|\bar{D}|$ are polynomially bounded when $K$ is fixed, as is the case for nanostructures ($K=3$).

    The tribond algorithm, however, requires {\it precise} distance data, while the distance data coming from experiments is noisy. The LIGA algorithm \cite{dg-4or} is a population-based heuristic designed to find realizations consistent with noisy distances or incomplete distance lists: it evaluates the fitness of each individual (partial) realization $x$ by the {\it mean square distance error}, i.e.
    \[\min_{\alpha} \frac{1}{m} \sum\limits_{\{u,v\}\in E} (\|x_u-x_v\|_2 - d_{\alpha^{-1}(u,v)})^2.\]

Given its practical importance, we believe this problem requires more work. Specifically, given that we possess theoretically and practically efficient methods for solving assignment problems \cite{assignment} and some practically efficient methods for solving the {\it i}DGP \cite{bpinterval}, the approach consisting in decomposing these two subproblems has not been sufficiently explored.

\subsubsection{Protein shape from NOESY data}
\label{s:noesy}
The NMR experiments that allow the estimation of inter-atomic distances in proteins are of the Nuclear Overhauser Effect (NOE) type: they are collectively called ``NOE spectrometry'', or NOESY for short. The actual NOESY output looks like a two dimensional surface with some peaks at some positions on the plane, which has axes labelled by chemical shift, a relative frequency measured in part-per-million (ppm). The peak intensity is related with the distance value arising in atoms within atomic groups that resonate at the given ppm values. 

We formalize this problem as follows \cite{2dnoesy}. Let $V$ be a set of atoms, and $\mathcal{I}=\{I_p\subset V\;|\;p<r\}$ be a given system of $r$ subsets of $V$ (representing the atomic groups resonating at a given chemical shift), so that
\[ \bigcup\limits_{p<r} I_p \subseteq V. \]
In general, $\mathcal{I}$ may not be a cover, i.e.~the above containment relation might be strict. Let $F=\{\{u,v\}\;|\;\exists p,q\le r\;(u\in I_p\land v\in I_q)\}$ be all of the possible edges between atoms in the atomic groups. Let $B$ be the set of pairs $\{u,v\}$ of atoms where the distances $d_{uv}$ are known (fairly) precisely --- because, e.g., they arise from covalent bonds and covalent angles.

The set system $\mathcal{I}$ is a proxy to the chemical shift value that index the axes of the NOE spectrum density distribution. The densitiy peak at coordinates $(I_p,I_q)$ is contained in the interval $[\rho^L_{pq},\rho^U_{pq}]$, which is what the NOESY experiments measure (for all $p<q\le r$). Consider the system:
\begin{eqnarray}
  \forall p<q\le r \qquad \rho_{pq}^L \le \sum\limits_{u\in I_p\atop v\in I_q} \frac{\beta_{uv}}{d_{uv}^6} \le \rho_{pq}^U \label{eq:noesyA} \\
  \forall \{u,v\}\in F\cup B \qquad \|x_u-x_v\|_2=d_{uv},\label{eq:noesyB}
\end{eqnarray}
where $\beta_{uv}$ are given constants. We remark that, in Eq.~\eqref{eq:noesyA}-\eqref{eq:noesyB}, the symbols $d_{uv}$ are decision variables rather than given constants.

The constants $\beta_{uv}$ are usually computed as averages over time of functions of angles between $x_u$ and $x_v$, but since $\beta_{uv}$ all converge to zero exponentially fast, it suffices to consider a discretization of $[0,\beta^U_{uv}]$ to derive a sampling for $\beta$. Thus $\beta$ can be assumed to be a constant (given) parameter tensor.

Finding solutions $(x,d)$ to Eq.~\eqref{eq:noesyA}-\eqref{eq:noesyB} for each slice of the tensor $\beta$ is the problem behind the NOESY assignment. This assumes that the ``sequential assignment problem'' has already been solved --- in other words, it assumes that we know what atoms (type and multiplicity) constitute the molecule. As far as we know, this is the first time that the problem has been formalized in this way. This formalization was reached in collaboration with two bioinformaticians of the structural biology group at Institut Pasteur, Dr.~T.~Malliavin and Dr.~B.~Worley. 

This problem is an interesting hybrid between DGP and uDGP: each distance value arises in correspondence to two chemical shift values, which in turn indicate a group of atoms. Moreover, distance data is measured by intervals, rather than precisely. More precisely, the intervals have an inverse sixth-power relationship with the distance values themselves, as shown in Eq.~\eqref{eq:noesyA}. 

The problem of protein structure determination from 2D NOESY spectra has traditionally been solved using metaheuristics, specifically simulated annealing and variants thereof \cite{dekkers,locatelli2002,nilges3}. More modern approaches employ convex optimization and graph decomposition into rigid components \cite{cucuringudgp,singer3-ii}. Differently than for nanostructures (Sect.~\ref{s:nano}), however, all currently existing methods appear to consider the problem as already decomposed into well distinguished ``assignment'' and ``realization'' phases \cite{khoo,donald}. We believe that a fresh look at the integrated problem (e.g.~in the formulation Eq.~\eqref{eq:noesyA}-\eqref{eq:noesyB} might suggest new methods, or at least new ways of decomposing the problem. 

\subsubsection{Relevance}
The uDGP is relevant because of the importance of its applications. The determination of nanostructures is crucial for material sciences, specifically for new, synthetic materials with exceptional properties in strength, stiffness or even auxeticity \cite{auxetics}. Interesting, hybrid problems between DGP and uDGP arise in another fundamentally important problem, that of estimating the geometrical shape of proteins from information obtained in 2D NOESY experiments \cite{wuthrich}. 

\section{Conclusion}
We have surveyed the field of distance geometry from the point of view of some open research areas which we think are very important: the search for a purely combinatorial characterization of three-dimensional rigidity, the determination of the complexity status of the Euclidean distance matrix completion problem, an investigation on techniques for estimating the number and type of solutions of various distance geometry problems, and progress on solution techniques for unassigned distance geometry problems.

% end paper

\section*{Acknowledgments}
We are grateful for Ana Flavia Lima for helping to check the manuscript prior to submission.
LL was partly supported by the ANR ``Bip:Bip'' project n.~ANR-10-BINF-03-08. CL was partly supported by the Brazilian research agencies FAPESP, CNPq, CAPES. 

\bibliographystyle{plain}
\bibliography{dr2}

\begin{thebibliography}{100}

\bibitem{aaronson}
S.~Aaronson.
\newblock Is \textbf{P} versus \textbf{NP} formally independent?
\newblock {\em Bulletin of the EATCS}, 81:Computational Complexity Column,
  2003.

\bibitem{abbott}
T.~Abbott.
\newblock Generalizations of kempe's universality theorem.
\newblock Master's thesis, MIT, 2008.

\bibitem{alexandrov}
A.~Alexandrov.
\newblock {\em Convex Polyhedra (in {R}ussian)}.
\newblock Gosudarstv.~Izdat.~Tekhn.-Theor.~Lit., Moscow, 1950.

\bibitem{alfakihdgp}
A.~Alfakih.
\newblock Universal rigidity of bar frameworks in general position: a
  {E}uclidean {D}istance {M}atrix approach.
\newblock In Mucherino et~al. \cite{dgpbook}, pages 3--22.

\bibitem{ipmsdp}
F.~Alizadeh.
\newblock Interior point methods in semidefinite programming with applications
  to combinatorial optimization.
\newblock {\em SIAM Journal on Optimization}, 5(1):13--51, 1995.

\bibitem{allender}
E.~Allender, R.~Beals, and M.~Ogihara.
\newblock The complexity of matrix rank and feasible systems of linear
  equations.
\newblock {\em Computational Complexity}, 8:99--126, 1999.

\bibitem{clifford-dga13}
R.~Alves, A.~Cassioli, A.~Mucherino, C.~Lavor, and L.~Liberti.
\newblock Adaptive branching in {\it i}{BP} with {C}lifford algebra.
\newblock In A.~Andrioni, C.~Lavor, L.~Liberti, A.~Mucherino, N.~Maculan, and
  R.~Rodriguez, editors, {\em Proceedings of the workshop on Distance Geometry
  and Applications}, Manaus, 2013. Universidade Federal do Amazonas.

\bibitem{asimow1}
L.~Asimow and B.~Roth.
\newblock The rigidity of graphs.
\newblock {\em Transactions of the American Mathematical Society},
  245:279--289, 1978.

\bibitem{asimow2}
L.~Asimow and B.~Roth.
\newblock The rigidity of graphs {II}.
\newblock {\em Journal of Mathematical Analysis and Applications}, 68:171--190,
  1979.

\bibitem{babai}
L.~Babai.
\newblock Automorphism groups, isomorphism, reconstruction.
\newblock In R.~Graham, M.~Gr\"otschel, and L.~Lov\'asz, editors, {\em Handbook
  of Combinatorics, vol.~2}, pages 1447--1540. MIT Press, Cambridge, MA, 1996.

\bibitem{babai2}
L.~Babai.
\newblock Graph isomorphism in quasipolynomial time.
\newblock Technical Report 1512.03547v2, ar{X}iv, 2016.

\bibitem{bahr}
A.~Bahr, J.~Leonard, and M.~Fallon.
\newblock Cooperative localization for autonomous underwater vehicles.
\newblock {\em International Journal of Robotics Research}, 28(6):714--728,
  2009.

\bibitem{barvinok}
A.~Barvinok.
\newblock Problems of distance geometry and convex properties of quadratic
  maps.
\newblock {\em Discrete and Computational Geometry}, 13:189--202, 1995.

\bibitem{barvinok2}
A.~Barvinok.
\newblock Measure concentration in optimization.
\newblock {\em Mathematical Programming}, 79:33--53, 1997.

\bibitem{pollack}
S.~Basu, R.~Pollack, and M.-F. Roy.
\newblock {\em Algorithms in real algebraic geometry}.
\newblock Springer, New York, 2006.

\bibitem{dgpinnp}
N.~Beeker, S.~Gaubert, C.~Glusa, and L.~Liberti.
\newblock Is the distance geometry problem in {{\bf NP}}?
\newblock In Mucherino et~al. \cite{dgpbook}, pages 85--94.

\bibitem{benedetti}
R.~Benedetti and J.-J. Risler.
\newblock {\em Real algebraic and semi-algebraic sets}.
\newblock Hermann, Paris, 1990.

\bibitem{berger}
B.~Berger, J.~Kleinberg, and T.~Leighton.
\newblock Reconstructing a three-dimensional model with arbitrary errors.
\newblock {\em Journal of the ACM}, 46(2):212--235, 1999.

\bibitem{dg-4or}
S.~Billinge, P.~Duxbury, D.~Gon\c{c}alves, C.~Lavor, and A.~Mucherino.
\newblock Assigned and unassigned distance geometry: {A}pplications to
  biological molecules and nanostructures.
\newblock {\em 4OR}, to appear.

\bibitem{biswas2004}
P.~Biswas and Y.~Ye.
\newblock Semidefinite programming for ad hoc wireless sensor network
  localization.
\newblock In {\em Proceedings of the 3rd international symposium on Information
  processing in sensor networks (IPSN04)}, pages 46--54, New York, NY, USA,
  2004. ACM.

\bibitem{blum}
L.~Blum, M.~Shub, and S.~Smale.
\newblock On a theory of computation and complexity over the real numbers: {\it
  NP}-completeness, recursive functions, and universal machines.
\newblock {\em Bulletin of the AMS}, 21(1):1--46, 1989.

\bibitem{blumenthal61}
L.~Blumenthal.
\newblock {\em A modern view of geometry}.
\newblock Freeman \& C., San Francisco, 1961.

\bibitem{bollobasrndgph}
B.~Bollob\'as.
\newblock {\em Random Graphs}.
\newblock Number~73 in Cambridge Studies in Advanced Mathematics. Cambridge
  University Press, Cambridge, 2001.

\bibitem{auxetics}
C.~Borcea and I.~Streinu.
\newblock Geometric auxetics.
\newblock {\em Proceedings of the Royal Society A}, 471(2184):20150033, 2015.

\bibitem{borg_10}
I.~Borg and P.~Groenen.
\newblock {\em Modern Multidimensional Scaling}.
\newblock Springer, New York, second edition, 2010.

\bibitem{bourgain}
J.~Bourgain.
\newblock On {L}ipschitz embeddings of finite metric spaces in {H}ilbert space.
\newblock {\em Israel Journal of Mathematics}, 52(1-2):46--52, 1985.

\bibitem{bowers}
J.C. Bowers and P.~Bowers.
\newblock A menger redux: Embedding metric spaces isometrically.
\newblock {\em American Mathematical Monthly}, in revision.

\bibitem{burgisser}
P.~B\"urgisser, M.~Clausen, and M.~Shokrollahi.
\newblock {\em Algebraic Complexity Theory}.
\newblock Number 315 in Grundlehren der mathematischen Wissenschaften.
  Springer, Berlin, 1997.

\bibitem{assignment}
R.~Burkard, M.~Dell'Amico, and S.~Martello.
\newblock {\em Assignment problems}.
\newblock SIAM, Providence, 2009.

\bibitem{bipbip}
A.~Cassioli, B.~Bordeaux, G.~Bouvier, A.~Mucherino, R.~Alves, L.~Liberti,
  M.~Nilges, C.~Lavor, and T.~Malliavin.
\newblock An algorithm to enumerate all possible protein conformations
  verifying a set of distance constraints.
\newblock {\em BMC Bioinformatics}, page 16:23, 2015.

\bibitem{cauchyrigid}
A.-L. Cauchy.
\newblock Sur les polygones et les poly\`edres.
\newblock {\em Journal de l'\'Ecole Polytechnique}, 16(9):87--99, 1813.

\bibitem{cayley1841}
A.~Cayley.
\newblock A theorem in the geometry of position.
\newblock {\em Cambridge Mathematical Journal}, II:267--271, 1841.

\bibitem{cheung}
H.Y. Cheung, T.C. Kwok, and L.C. Lau.
\newblock Fast matrix rank algorithms and applications.
\newblock {\em Journal of the ACM}, 60(5):31:1--31:25, 2013.

\bibitem{cobham}
A.~Cobham.
\newblock The intrinsic computational difficulty of functions.
\newblock In Y.~Bar-Hillel, editor, {\em Logic, Methodology and Philosophy of
  Science}, pages 24--30. North-Holland, Amsterdam, 1965.

\bibitem{connelly-countereg}
R.~Connelly.
\newblock A counterexample to the rigidity conjecture for polyhedra.
\newblock {\em Publications Math\'ematiques de l'IHES}, 47:333--338, 1978.

\bibitem{cook}
S.~Cook.
\newblock The complexity of theorem-proving procedures.
\newblock In {\em ACM Symposium on the Theory of Computing}, STOC, pages
  151--158, New York, 1971. ACM.

\bibitem{coope}
I.~Coope.
\newblock Reliable computation of the points of intersection of $n$ spheres in
  $\mathbb{R}^n$.
\newblock {\em Australian and New Zealand Industrial and Applied Mathematics
  Journal}, 42:C461--C477, 2000.

\bibitem{sidechains}
V.~Costa, A.~Mucherino, C.~Lavor, A.~Cassioli, L.~Carvalho, and N.~Maculan.
\newblock Discretization orders for protein side chains.
\newblock {\em Journal of Global Optimization}, 60:333--349, 2014.

\bibitem{cremona1872}
L.~Cremona.
\newblock {\em Le figure reciproche nella statica grafica}.
\newblock G.~Bernardoni, Milano, 1872.

\bibitem{cremona1874}
L.~Cremona.
\newblock {\em Elementi di calcolo grafico}.
\newblock Paravia, Torino, 1874.

\bibitem{crippeninf}
G.~Crippen.
\newblock An alternative approach to distance geometry using $l^\infty$
  distances.
\newblock {\em Discrete Applied Mathematics}, 197:20--26, 2015.

\bibitem{cucuringudgp}
M.~Cucuringu.
\newblock Asap -- an eigenvector synchronization algorithm for the graph
  realization problem.
\newblock In Mucherino et~al. \cite{dgpbook}, pages 177--196.

\bibitem{singer3-ii}
M.~Cucuringu, A.~Singer, and D.~Cowburn.
\newblock Eigenvector synchronization, graph ridigity and the molecule problem.
\newblock {\em Information and Inference: a journal of the {IMA}}, 1:21--67,
  2012.

\bibitem{dakic}
T.~Daki\'c.
\newblock {\em On the turnpike problem}.
\newblock PhD thesis, Simon Fraser University, 2000.

\bibitem{dgpzoo-tr}
C.~D'Ambrosio, Vu~Ky, C.~Lavor, L.~Liberti, and N.~Maculan.
\newblock Solving distance geometry problems with interval data using
  formulation-based methods.
\newblock Technical report, LIX Ecole Polytechnique (working paper), 2014.

\bibitem{oneinfnorm}
C.~D'Ambrosio and L.~Liberti.
\newblock Distance geometry in linearizable norms.
\newblock Technical Report working paper, Ecole Polytechnique, 2016.

\bibitem{dekkers}
A.~Dekkers and E.~Aarts.
\newblock Global optimization and simulated annealing.
\newblock {\em Mathematical Programming}, 50:367--393, 1991.

\bibitem{isco16}
G.~Dias and L.~Liberti.
\newblock Diagonally dominant programming in distance geometry.
\newblock In R.~Cerulli, S.~Fujishige, and R.~Mahjoub, editors, {\em
  International Symposium in Combinatorial Optimization}, volume 9849 of {\em
  LNCS}, pages 225--236, New York, 2016. Springer.

\bibitem{doherty}
L.~Doherty, K.~Pister, and L.~{El Ghaoui}.
\newblock Convex position estimation in wireless sensor networks.
\newblock In {\em Twentieth Annual Joint Conference of the IEEE Computer and
  Communications Societies}, volume~3 of {\em INFOCOM}, pages 1655--1663,
  Piscataway, 2001. IEEE.

\bibitem{vetterli}
I.~Dokmani\'c, R.~Parhizkar, J.~Ranieri, and M.~Vetterli.
\newblock Euclidean distance matrices: Essential theory, algorithms and
  applications.
\newblock {\em IEEE Signal Processing Magazine}, 1053-5888:12--30, Nov. 2015.

\bibitem{donald}
B.~Donald.
\newblock {\em Algorithms in Structural Molecular Biology}.
\newblock MIT Press, Boston, 2011.

\bibitem{dongwu}
Q.~Dong and Z.~Wu.
\newblock A linear-time algorithm for solving the molecular distance geometry
  problem with exact inter-atomic distances.
\newblock {\em Journal of Global Optimization}, 22:365--375, 2002.

\bibitem{splogic}
H.~Du, N.~Alechina, K.~Stock, and M.~Jackson.
\newblock The logic of {NEAR} and {FAR}.
\newblock In T.~Tenbrink et~al., editor, {\em COSIT}, volume 8116 of {\em
  LNCS}, pages 475--494, Switzerland, 2013. Springer.

\bibitem{tribond}
P.~Duxbury, L.~Granlund, P.~Juhas, and S.~Billinge.
\newblock The unassigned distance geometry problem.
\newblock {\em Discrete Applied Mathematics}, 204:117--132, 2016.

\bibitem{edmonds}
J.~Edmonds.
\newblock Paths, trees and flowers.
\newblock {\em Canadian Journal of Mathematics}, 17:449--467, 1965.

\bibitem{erdos-renyi}
P.~Erd\H{o}s and A.~Renyi.
\newblock On the evolution of random graphs.
\newblock {\em Publications of the Mathematical Institute of the Hungarian
  Academy of Sciences}, 5:17--61, 1960.

\bibitem{eren04}
T.~Eren, D.~Goldenberg, W.~Whiteley, Y.~Yang, A.~Morse, B.~Anderson, and
  P.~Belhumeur.
\newblock Rigidity, computation, and randomization in network localization.
\newblock {\em IEEE}, pages 2673--2684, 2004.

\bibitem{euler1766}
L.~Euler.
\newblock Continuatio fragmentorum ex adversariis mathematicis depromptorum:
  {II} {G}eometria, 97.
\newblock In P.~Fuss and N.~Fuss, editors, {\em Opera postuma mathematica et
  physica anno 1844 detecta}, volume~I, pages 494--496. Eggers \& C.,
  Petropolis, 1862.

\bibitem{encopt2}
C.~Floudas and P.~Pardalos, editors.
\newblock {\em Encyclopedia of Optimization}.
\newblock Springer, New York, second edition, 2009.

\bibitem{frechet}
M.~Fr\'echet.
\newblock Sur quelques points du calcul fonctionnel.
\newblock {\em Rendiconti del Circolo Matematico di Palermo}, 22:1--74, 1906.

\bibitem{gareyjohnson}
M.~Garey and D.~Johnson.
\newblock {\em Computers and Intractability: a Guide to the Theory of {{\bf
  NP}}-{C}ompleteness}.
\newblock Freeman and Company, New York, 1979.

\bibitem{garibaldi}
J.~Garibaldi, A.~Iosevich, and S.~Senger.
\newblock {\em The Erd\H{o}s Distance Problem}.
\newblock Number~56 in Student Mathematical Library. AMS, Providence, 2011.

\bibitem{gluck}
H.~Gluck.
\newblock Almost all simply connected closed surfaces are rigid.
\newblock In A.~Dold and B.~Eckmann, editors, {\em Geometric Topology}, volume
  438 of {\em Lecture Notes in Mathematics}, pages 225--239, Berlin, 1975.
  Springer.

\bibitem{goedelDG1}
K.~G\"odel.
\newblock On the isometric embeddability of quadruples of points of $r_3$ in
  the surface of a sphere.
\newblock In S.~Feferman, J.~Dawson, S.~Kleene, G.~Moore, R.~Solovay, and
  J.~van Heijenoort, editors, {\em Kurt G\"odel: Collected Works, vol.~I},
  pages (1933b) 276--279. Oxford University Press, Oxford, 1986.

\bibitem{gortler}
S.~Gortler, A.~Healy, and D.~Thurston.
\newblock Characterizing generic global rigidity.
\newblock {\em American Journal of Mathematics}, 132(4):897--939, 2010.

\bibitem{graverbook}
J.~Graver, B.~Servatius, and H.~Servatius.
\newblock {\em Combinatorial Rigidity}.
\newblock American Mathematical Society, 1993.

\bibitem{grunbaum}
B.~Gr\"unbaum and G.~Shephard.
\newblock Lectures on lost mathematics.
\newblock Technical Report EPrint Collection -- Mathematics [112], University
  of Washington, 2010.

\bibitem{havel}
T.~Havel, I.~Kuntz, and G.~Crippen.
\newblock The theory and practice of distance geometry.
\newblock {\em Bulletin of Mathematical Biology}, 45(5):665--720, 1983.

\bibitem{wuthrich}
T.~Havel and K.~W\"uthrich.
\newblock An evaluation of the combined use of nuclear magnetic resonance and
  distance geometry for the determination of protein conformations in solution.
\newblock {\em Journal of Molecular Biology}, 182(2):281--294, 1985.

\bibitem{Hen92}
B.~Hendrickson.
\newblock Conditions for unique graph realizations.
\newblock {\em SIAM Journal on Computing}, 21(1):65--84, 1992.

\bibitem{heron}
Heron.
\newblock {\em Metrica}, volume~I.
\newblock Alexandria, $\sim$\!100AD.

\bibitem{hoang}
T.M. Hoang.
\newblock {\em On the complexity of some problems in linear algebra}.
\newblock PhD thesis, Universit\"at Ulm, 2003.

\bibitem{indyk}
P.~Indyk and A.~Naor.
\newblock Nearest neighbor preserving embeddings.
\newblock {\em ACM Transactions on Algorithms}, 3(3):Art.~31, 2007.

\bibitem{globrigid2}
B.~Jackson and T.~Jord\'an.
\newblock Connected rigidity matroids and unique realization of graphs.
\newblock {\em Journal of Combinatorial Theory, Series B}, 94:1--29, 2005.

\bibitem{jjsurvey}
B.~Jackson and T.~Jord\'an.
\newblock Graph theoretic techniques in the analysis of uniquely localizable
  sensor networks.
\newblock In G.~Mao and B.~Fidan, editors, {\em Localization Algorithms and
  Strategies for Wireless Sensor Networks: Monitoring and Surveillance
  Techniques for Target Tracking}, pages 146--173. IGI Global, 2009.

\bibitem{thorpe_perc}
D.~Jacobs and M.~Thorpe.
\newblock Generic rigidity percolation.
\newblock {\em Physical Review Letters}, 75(22):4051--4054, 1995.

\bibitem{johnson1982}
D.~Johnson.
\newblock The {{\bf NP}}-completeness column: an ongoing guide.
\newblock {\em Journal of Algorithms}, 3:182--195, 1982.

\bibitem{jllemma}
W.~Johnson and J.~Lindenstrauss.
\newblock Extensions of {L}ipschitz mappings into a {H}ilbert space.
\newblock In G.~Hedlund, editor, {\em Conference in Modern Analysis and
  Probability}, volume~26 of {\em Contemporary Mathematics}, pages 189--206,
  Providence, 1984. American Mathematical Society.

\bibitem{karp}
R.~Karp.
\newblock Reducibility among combinatorial problems.
\newblock In R.~Miller and W.~Thatcher, editors, {\em Complexity of Computer
  Computations}, volume~5 of {\em IBM Research Symposia}, pages 85--104, New
  York, 1972. Plenum.

\bibitem{khoo}
Yuehaw Khoo.
\newblock {\em Protein structural calculation from {NMR} spectroscopy}.
\newblock PhD thesis, Princeton University, 2016.

\bibitem{krislocksiam}
N.~Krislock and H.~Wolkowicz.
\newblock Explicit sensor network localization using semidefinite
  representations and facial reductions.
\newblock {\em SIAM Journal on Optimization}, 20:2679--2708, 2010.

\bibitem{kuratowski}
C.~Kuratowski.
\newblock Quelques probl\`emes concernant les espaces m\'etriques
  non-s\'eparables.
\newblock {\em Fundamenta Mathematic{\ae}}, 25:534--545, 1935.

\bibitem{laman}
G.~Laman.
\newblock On graphs and rigidity of plane skeletal structures.
\newblock {\em Journal of Engineering Mathematics}, 4(4):331--340, 1970.

\bibitem{laurent00}
M.~Laurent.
\newblock Polynomial instances of the positive semidefinite and {E}uclidean
  distance matrix completion problems.
\newblock {\em SIAM Journal of Matrix Analysis and Applications},
  22(3):874--894, 2000.

\bibitem{mcp}
M.~Laurent.
\newblock Matrix completion problems.
\newblock In Floudas and Pardalos \cite{encopt2}, pages 1967--1975.

\bibitem{Lav05}
C.~Lavor.
\newblock On generating instances for the molecular distance geometry problem.
\newblock In L.~Liberti and N.~Maculan, editors, {\em Global Optimization: from
  Theory to Implementation}, pages 405--414. Springer, Berlin, 2006.

\bibitem{cliffordalgebra}
C.~Lavor, R.~Alves, W.~Figuereido, A.~Petraglia, and N.~Maculan.
\newblock Clifford algebra and the discretizable molecular distance geometry
  problem.
\newblock {\em Advances in Applied Clifford Algebras}, 25:925--942, 2015.

\bibitem{dgpitorpreface}
C.~Lavor, M.~Firer, J.-M. Martinez, and L.~Liberti.
\newblock Preface.
\newblock {\em International Transactions in Operational Research}, 23(5):841,
  2016.

\bibitem{dvop}
C.~Lavor, J.~Lee, A.~{Lee-St.~John}, L.~Liberti, A.~Mucherino, and
  M.~Sviridenko.
\newblock Discretization orders for distance geometry problems.
\newblock {\em Optimization Letters}, 6:783--796, 2012.

\bibitem{lln4}
C.~Lavor, L.~Liberti, and N.~Maculan.
\newblock Molecular distance geometry problem.
\newblock In Floudas and Pardalos \cite{encopt2}, pages 2305--2311.

\bibitem{dmdgp}
C.~Lavor, L.~Liberti, N.~Maculan, and A.~Mucherino.
\newblock The discretizable molecular distance geometry problem.
\newblock {\em Computational Optimization and Applications}, 52:115--146, 2012.

\bibitem{dmdgpejor}
C.~Lavor, L.~Liberti, N.~Maculan, and A.~Mucherino.
\newblock Recent advances on the discretizable molecular distance geometry
  problem.
\newblock {\em European Journal of Operational Research}, 219:698--706, 2012.

\bibitem{bpinterval}
C.~Lavor, L.~Liberti, and A.~Mucherino.
\newblock The {\it interval} {Branch-and-Prune} algorithm for the discretizable
  molecular distance geometry problem with inexact distances.
\newblock {\em Journal of Global Optimization}, 56:855--871, 2013.

\bibitem{skiena}
P.~Lemke, S.~Skiena, and W.~Smith.
\newblock Reconstructing sets from interpoint distances.
\newblock In B.~Aronov et~al., editor, {\em Discrete and Computational
  Geometry}, volume~25 of {\em Algorithms and Combinatorics}, pages 597--631,
  Berlin, 2003. Springer.

\bibitem{six}
L.~Liberti and C.~Lavor.
\newblock Six mathematical gems in the history of distance geometry.
\newblock {\em International Transactions in Operational Research},
  23:897--920, 2016.

\bibitem{liberti-gsi13}
L.~Liberti, C.~Lavor, J.~Alencar, and G.~Abud.
\newblock Counting the number of solutions of ${}^k${DMDGP} instances.
\newblock In F.~Nielsen and F.~Barbaresco, editors, {\em Geometric Science of
  Information}, volume 8085 of {\em LNCS}, pages 224--230, New York, 2013.
  Springer.

\bibitem{lln5}
L.~Liberti, C.~Lavor, and N.~Maculan.
\newblock A branch-and-prune algorithm for the molecular distance geometry
  problem.
\newblock {\em International Transactions in Operational Research}, 15:1--17,
  2008.

\bibitem{dgp-sirev}
L.~Liberti, C.~Lavor, N.~Maculan, and A.~Mucherino.
\newblock Euclidean distance geometry and applications.
\newblock {\em SIAM Review}, 56(1):3--69, 2014.

\bibitem{bppolybook}
L.~Liberti, C.~Lavor, and A.~Mucherino.
\newblock The discretizable molecular distance geometry problem seems easier on
  proteins.
\newblock In Mucherino et~al. \cite{dgpbook}, pages 47--60.

\bibitem{mdgpsurvey}
L.~Liberti, C.~Lavor, A.~Mucherino, and N.~Maculan.
\newblock Molecular distance geometry methods: from continuous to discrete.
\newblock {\em International Transactions in Operational Research}, 18:33--51,
  2010.

\bibitem{powerof2}
L.~Liberti, B.~Masson, C.~Lavor, J.~Lee, and A.~Mucherino.
\newblock On the number of realizations of certain {H}enneberg graphs arising
  in protein conformation.
\newblock {\em Discrete Applied Mathematics}, 165:213--232, 2014.

\bibitem{dgpsphere}
L.~Liberti, G.~Swirszcz, and C.~Lavor.
\newblock Distance geometry on the sphere.
\newblock In {\em JCDCG${}^2$}, LNCS, New York, accepted. Springer.

\bibitem{dgta16proc}
Leo Liberti, editor.
\newblock {\em Proceedings of the DIMACS Workshop on Distance Geometry Theory
  and Applications (DGTA16)}, 2016.

\bibitem{linial-stoc}
N.~Linial, E.~London, and Y.~Rabinovich.
\newblock The geometry of graphs and some of its algorithmic applications.
\newblock In {\em Proceedings of the Symposium on Foundations of Computer
  Science}, volume~35 of {\em FOCS}, pages 577--591, Piscataway, 1994. IEEE.

\bibitem{linial}
N.~Linial, E.~London, and Y.~Rabinovich.
\newblock The geometry of graphs and some of its algorithmic applications.
\newblock {\em Combinatorica}, 15(2):215--245, 1995.

\bibitem{locatelli2002}
M.~Locatelli.
\newblock Simulated annealing algorithms for global optimization.
\newblock In P.M. Pardalos and H.E. Romeijn, editors, {\em Handbook of Global
  Optimization}, volume~2, pages 179--229. Kluwer Academic Publishers,
  Dordrecht, 2002.

\bibitem{lovasz-yemini}
L.~Lov\'asz and Y.~Yemini.
\newblock On generic rigidity in the plane.
\newblock {\em SIAM Journal on Algebraic and Discrete Methods}, 3(1):91--98,
  1982.

\bibitem{mahajansarma}
M.~Mahajan and J.~Sarma.
\newblock On the complexity of matrix rank and rigidity.
\newblock {\em Theory of Computing Systems}, 46:9--26, 2010.

\bibitem{2dnoesy}
T.~Malliavin and B.~Worley.
\newblock Personal communication, 2016.
\newblock Institut Pasteur, Paris, France.

\bibitem{matousekmetric}
J.~Matou\v{s}ek.
\newblock Lecture notes on metric embeddings.
\newblock Technical report, ETH Z\"urich, 2013.

\bibitem{maxwell1864b}
J.~Maxwell.
\newblock On reciprocal figures and diagrams of forces.
\newblock {\em Philosophical Magazine}, 27(182):250--261, 1864.

\bibitem{maxwell1864}
J.~Maxwell.
\newblock On the calculation of the equilibrium and stiffness of frames.
\newblock {\em Philosophical Magazine}, 27(182):294--299, 1864.

\bibitem{mehlhorn}
K.~Mehlhorn and P.~Sanders.
\newblock {\em Algorithms and Data Structures}.
\newblock Springer, Berlin, 2008.

\bibitem{menger28}
K.~Menger.
\newblock {U}ntersuchungen \"uber allgemeine {M}etrik.
\newblock {\em Mathematische Annalen}, 100:75--163, 1928.

\bibitem{menger31}
K.~Menger.
\newblock New foundation of {E}uclidean geometry.
\newblock {\em American Journal of Mathematics}, 53(4):721--745, 1931.

\bibitem{milnor64}
J.~Milnor.
\newblock On the {B}etti numbers of real varieties.
\newblock {\em Proceedings of the American Mathematical Society}, 15:275--280,
  1964.

\bibitem{morewu}
J.~Mor\'e and Z.~Wu.
\newblock Global continuation for distance geometry problems.
\newblock {\em SIAM Journal of Optimization}, 7(3):814--846, 1997.

\bibitem{mosek7}
Mosek ApS.
\newblock {\em The {\tt mosek} manual, Version 7 (Revision 114)}, 2014.
\newblock \verb+(www.mosek.com)+.

\bibitem{duxbury_perc}
C.~Moukarzel and P.~Duxbury.
\newblock Stressed backbone and elasticity of random central-force systems.
\newblock {\em Physical Review Letters}, 75(22):4055--4059, 1995.

\bibitem{dgpdampreface}
A.~Mucherino, R.~de~Freitas, and C.~Lavor.
\newblock Preface.
\newblock {\em Discrete Applied Mathematics}, 197:1--2, 2015.

\bibitem{ddgp}
A.~Mucherino, C.~Lavor, and L.~Liberti.
\newblock The discretizable distance geometry problem.
\newblock {\em Optimization Letters}, 6:1671--1686, 2012.

\bibitem{symmBPjbcb}
A.~Mucherino, C.~Lavor, and L.~Liberti.
\newblock Exploiting symmetry properties of the discretizable molecular
  distance geometry problem.
\newblock {\em Journal of Bioinformatics and Computational Biology},
  10:1242009(1--15), 2012.

\bibitem{dgpbook}
A.~Mucherino, C.~Lavor, L.~Liberti, and N.~Maculan, editors.
\newblock {\em Distance Geometry: Theory, Methods, and Applications}.
\newblock Springer, New York, 2013.

\bibitem{mdjeep}
A.~Mucherino, L.~Liberti, and C.~Lavor.
\newblock {\tt MD-jeep}: an implementation of a branch-and-prune algorithm for
  distance geometry problems.
\newblock In K.~Fukuda, J.~van~der Hoeven, M.~Joswig, and N.~Takayama, editors,
  {\em Mathematical Software}, volume 6327 of {\em LNCS}, pages 186--197, New
  York, 2010. Springer.

\bibitem{nilges3}
M.~Nilges, G.M. Clore, and A.~Gronenborn.
\newblock Determination of three-dimensional structures of proteins from
  interproton distance data by hybrid distance geometry-dynamical simulated
  annealing calculations.
\newblock {\em FEBS Letters}, 229(2):317--324, 1988.

\bibitem{nilges}
M.~Nilges, A.M. Gronenborn, A.T. Brunger, and G.M. Clore.
\newblock Determination of three-dimensional structures of proteins by
  simulated annealing with interproton distance restraints. application to
  crambin, potato carboxypeptidase inhibitor and barley serine proteinase
  inhibitor 2.
\newblock {\em Protein Engineering}, 2:27--38, 1988.

\bibitem{recskisurvey1}
A.~Recski.
\newblock Applications of combinatorics to statics --- {A} survey.
\newblock {\em Rendiconti del Circolo Matematico di Palermo},
  II(Suppl.~3):237--247, 1984.

\bibitem{recskisurvey2}
A.~Recski.
\newblock Applications of combinatorics to statics --- {A} second survey.
\newblock {\em Discrete Mathematics}, 108:183--188, 1992.

\bibitem{rojas_10}
N.~Rojas.
\newblock {\em Distance-based formulations for the position analysis of
  kinematic chains}.
\newblock PhD thesis, Universitat Politecnica de Catalunya, 2012.

\bibitem{santana07}
R.~Santana, P.~{Larra{\~{n}}aga}, and J.~Lozano.
\newblock Side chain placement using estimation of distribution algorithms.
\newblock {\em Artificial Intelligence in Medicine}, 39:49--63, 2007.

\bibitem{santana08}
R.~Santana, P.~{Larra{\~{n}}aga}, and J.A. Lozano.
\newblock Combining variable neighbourhood search and estimation of
  distribution algorithms in the protein side chain placement problem.
\newblock {\em Journal of Heuristics}, 14:519--547, 2008.

\bibitem{saxe79}
J.~Saxe.
\newblock Embeddability of weighted graphs in $k$-space is strongly {{\bf
  NP}}-hard.
\newblock {\em Proceedings of 17th Allerton Conference in Communications,
  Control and Computing}, pages 480--489, 1979.

\bibitem{schoenberg}
I.~Schoenberg.
\newblock Remarks to {M}aurice {F}r\'echet's article ``{S}ur la d\'efinition
  axiomatique d'une classe d'espaces distanci\'es vectoriellement applicable
  sur l'espace de {H}ilbert".
\newblock {\em Annals of Mathematics}, 36(3):724--732, 1935.

\bibitem{singer4}
A.~Singer.
\newblock Angular synchronization by eigenvectors and semidefinite programming.
\newblock {\em Applied and Computational Harmonic Analysis}, 30:20--36, 2011.

\bibitem{sitharam}
M.~Sitharam and Y.~Zhou.
\newblock A tractable, approximate, combinatorial 3{D} rigidity
  characterization.
\newblock In {\em Fifth Workshop on Automated Deduction in Geometry}, 2004.

\bibitem{sljoka}
A.~Sljoka.
\newblock {\em Algorithms in rigidity theory with applications to protein
  flexibility and mechanical linkages}.
\newblock PhD thesis, York University, Canada, 2012.

\bibitem{souza2}
M.~Souza, C.~Lavor, A.~Muritiba, and N.~Maculan.
\newblock Solving the molecular distance geometry problem with inaccurate
  distance data.
\newblock {\em BMC Bioinformatics}, 14(Suppl.~9):S71--S76, 2013.

\bibitem{tarski-reals}
A.~Tarski.
\newblock A decision method for elementary algebra and geometry.
\newblock Technical Report R-109, Rand Corporation, 1951.

\bibitem{tay-whiteley}
T.-S. Tay and W.~Whiteley.
\newblock Generating isostatic frameworks.
\newblock {\em Structural Topology}, 11:21--69, 1985.

\bibitem{thorpe}
M.~Thorpe and P.~Duxbury, editors.
\newblock {\em Rigidity Theory and Applications}.
\newblock Fundamental Materials Research. Springer, New York, 2002.

\bibitem{wuthrich_89}
K.~W\"uthrich.
\newblock Protein structure determination in solution by nuclear magnetic
  resonance spectroscopy.
\newblock {\em Science}, 243:45--50, 1989.

\bibitem{wuthrich_83}
K.~W\"uthrich, M.~Billeter, and W.~Braun.
\newblock Pseudo-structures for the 20 common amino acids for use in studies of
  protein conformations by measurements of intramolecular proton-proton
  distance constraints with nuclear magnetic resonance.
\newblock {\em Journal of Molecular Biology}, 169:949--961, 1983.

\bibitem{yemini78}
Y.~Yemini.
\newblock The positioning problem --- a draft of an intermediate summary.
\newblock In {\em Proceedings of the Conference on Distributed Sensor
  Networks}, pages 137--145, Pittsburgh, 1978. Carnegie-Mellon University.

\bibitem{yemini}
Y.~Yemini.
\newblock Some theoretical aspects of position-location problems.
\newblock In {\em Proceedings of the 20th Annual Symposium on the Foundations
  of Computer Science}, pages 1--8, Piscataway, 1979. IEEE.

\bibitem{young_householder}
G.~Young and A.~Householder.
\newblock Discussion of a set of points in terms of their mutual distances.
\newblock {\em Psychometrika}, 3(1):19--22, 1938.

\end{thebibliography}

\end{document}